\documentclass{amsart}
\usepackage{amsmath}
\usepackage{amssymb}
\catcode `\@=11
\def\numberbysection{\@addtoreset{equation}{section}
         \renewcommand{\theequation}{\thesection.\arabic{equation}}}
\numberbysection
\def\subsubsection{\@startsection{subsubsection}{3}%
  \normalparindent{.5\linespacing\@plus.7\linespacing}{-.5em}%    
  {\normalfont\bfseries}}
\def\cal{\mathcal}

\def\del{\partial}
\def\iso{\simeq}
\def\a{\alpha}

\def\g{\gamma}
\def\d{\delta}
\def\eps{\epsilon}
\def\l{\lambda}
\def\s{\sigma}
\def\t{\tau}
\def\Cal{\cal}  

\def\Zz{{\bf Z}/2{\bf Z}}

\def\Znn{{\bf Z}/(n+1){\bf Z}}

\def\nn{\nonumber}
\def\la{\langle}
\def\ra{\rangle}
\def\Z2{{\bf Z}/2{\bf Z}}

\def\Sn{{\bf S}_n}
\def\Snn{{\bf S}_{n+1}}
\def\sign{\mathrm{sign}}
\def\mathbb{\bf}
\def\codim{\rm{codim}}
\def\Hom{\rm{Hom}}
\def\C{\mathcal{C}}
\def\F{\mathcal{F}}
\def\Frob{\mathcal{FROB}}
\def\Tr{\mathrm{Tr}}
\def\STr{\mathrm{STr}}
\def\Fix{\mathrm{Fix}}
\def\sign{sign}

\begin{document}

\title{Second quantized Frobenius algebras} 

\author[Ralph M. Kaufmann]{Ralph M. Kaufmann$^*$ \\ 
University of Southern California, Los Angeles, USA\\
 and IH\'ES, Bures--sur--Yvette, France}

\thanks{Partially supported by NSF grant \#0070681}
\email{kaufmann@math.usc.edu}

\address{University of Southern California, Los Angeles, USA and
IH\'ES, Bures--sur--Yvette, France}

\begin{abstract}
We show that given a Frobenius algebra there is a unique notion
of its second quantization, which is the sum over all symmetric group 
quotients of n--th tensor powers, where the quotients are given
by symmetric group twisted Frobenius algebras.
To this end, 
we consider the setting of Frobenius algebras given by
functors from geometric categories whose objects are endowed with 
geometric group actions and prove structural results, which in
turn yield a constructive realization in the case of n--th tensor powers and
the natural permutation action. 
We also show that naturally graded symmetric
group twisted Frobenius algebras have a unique algebra structure
already determined by their underlying additive data together with
a choice of super--grading.
Furthermore we discuss several notions of discrete torsion and
show that indeed a non--trivial discrete torsion leads to a non--trivial
super structure on the second quantization. 
\end{abstract}

\maketitle

\tableofcontents
\section*{Introduction}

In ``stringy'' geometry evaluating a functor from a geometric to a linear
category on a group quotient is generally a two step process.
The first is to evaluate the functor not only on the object, but to 
form the direct sum of the evaluations on all of the fixed point sets.
The new summands corresponding to group elements which are not the identity
are usually named twisted sectors.
The second step is to find a suitable group action on the twisted sectors and
take group invariants.

If the objects in the linear category also have an algebra structure there
is an additional step, i.e.\ to find a new algebra structure that is not
the diagonal one which is canonically present, but a group graded one.

If there is also a natural pairing such that the original functor have 
values in Frobenius algebras, then the result of the ``stringy'' extension of 
this functor  should have values in  G--twisted Frobenius algebras which 
were introduced for this purpose in [K2].

In particular, the question of importance is the step of finding the suitable 
multiplication. The theory of G--twisted Frobenius algebras is exactly tailored
to classify the possible multiplicative structures.

We address this matter in the present paper once more in the general case
of intersection Frobenius algebras and in the 
special case of symmetric group quotients which are naturally intersection
Frobenius algebras.

The class of intersection Frobenius algebras incorporates the fact that all
geometric construction of Frobenius algebras via functors from 
geometric categories with geometric group actions actually have a much richer
structure which can be used to provide further constraints on the nature
of the twisted multiplication.

We apply these general results to the case of symmetric group quotients
of powers of Frobenius algebras.

The main result here is that {\em 
there is a unique multiplicative structure that
makes the canonical extension of the n--th tensor power of a Frobenius 
algebra into a symmetric group twisted Frobenius algebra}.

This uniqueness has
to be understood up to a twist by discrete torsion which is always
possible and up to a super re--grading. The former is parametrized
by $Z^2(\Sn,k^*)$ and up to isomorphism by $H^2(\Sn,k^*)=\Zz$
and the latter is also a choice in $\Zz$ which renders everything
either purely even or super--graded.

This result should be read as the statement that there is a well defined 
notion of second quantized Frobenius algebra. 
Recall that in the spirit of 
[DMVV,D1] second quantization in an monoidal category with a notion 
of symmetric quotients is given by:

Second quantization of X = $\exp(X)$ = $\sum_n X^{\times n}/\Sn$,\\
where $\Sn$ is acting by permutations on the factors and the sum may
be formal or contain a bookkeeping variable (e.g.\ $q^n$).
From our result, we expect that one can also easily 
derive a definition of second quantized 
motives. 
All the objects are powers of the original object
and the morphisms are given by structural morphisms.
It would be interesting to explicitly see the multiplication in terms of
correspondences.

Furthermore we discuss several notions of discrete torsion and
show that indeed a non--trivial discrete torsion leads to a non--trivial
super structure on the second quantization. 

The paper is organized as follows. In \S1 
we review our definitions of [K2, K3] of $G$--twisted and
special $G$--twisted Frobenius algebras. In
the latter the multiplication and group action can 
be described by group cocycles and non--abelian 
group cocycles, respectively. Besides fixing and recalling
the notation and definitions, we add several useful practical 
Lemmas as well as a new description of the non--abelian cocycles in terms
of ordinary group one--cocycles with values in tori.
The second paragraph contains the functorial setup of the general 
question posed in the introduction, i.e.\ to identify 
the underlying additive data and the possible extensions of this data
by ``stringy'' product to the right type of group quotient algebra.

In \S2 we also introduce the notion of intersection categories, which 
reflect the geometrical setups with geometrical group actions which 
are used for the known construction of Frobenius algebras such as
cohomology, quantum cohomology, singularity theory, etc..
This setup is carried over to the Frobenius side in \S 3 
where we prove general results about the structure of the cocycles in 
the special $G$--twisted Frobenius algebra case. These results 
are also the key to understanding the second quantization.
Furthermore we introduce the notion of algebraic discrete torsion, which
generalizes the case of discrete torsion for Jacobian algebras of [K3]
and provides the discrete torsion that is linked to the super--structure
of second quantization.

In order to give a clearer view of the geometry involved in the second
quantization, it is useful to also consider the case of Jacobian Frobenius
algebras and their second quantization. The relevant notions of Jacobian
Frobenius algebras are recalled in \S4.

We then start our consideration of
$\Sn$--twisted Frobenius algebras. \S 5 contains general results
about these structures. The main results of this section are 
the classification of possible non--abelian group cocycles and 
the uniqueness (up to normalization) of ``stringy'' products given
a group grading compatible with the natural grading on $\Sn$.
Before applying these results to general symmetric powers, we 
work out all the details in the case of the n--th tensor power of 
a Frobenius algebra in \S 6 and also show the existence of
the natural $\Sn$--twisted Frobenius algebra based on the  
n--th tensor power.
Here we also recover the known discrete torsion 
corresponding to the non--trivial Schur multiplier.

Using the geometric insight of the previous paragraph we 
turn to the general case of the  n--th tensor power  of 
a Frobenius algebra in \S7 and show that there is a unique (up to a choice of 
parity for the group action) natural extension of   n--th tensor power
to a  $\Sn$--twisted Frobenius algebra, establishing the existence of
second quantized Frobenius algebras. There are two versions,
a purely even one a and  super symmetric one.  
Passing from one to the other can be viewed as turning
on a natural algebraic discrete torsion.
Lastly, we relate our results 
to the ones of [LS].

There are also two appendices. The first contains a key result
on the possible form of non--abelian $\Sn$ cocycles and the
second contains the detailed version of the proof of normalizability
of \S 5.

\section*{Acknowledgments}
I would like to thank the IH\'ES for its kind hospitality.
My visits in 2001 and 2002 mark the conceptual origin and 
the finishing phase of the paper. I also gratefully 
acknowledge the support from the NSF. 
It is  a pleasure to 
thank L.\ Borisov who sparked my interested in symmetric products, 
Y.\ Ruan and A.\ Adem for discussions and the  wonderful
conference in Madison and B.\ Guralnick for discussions on 
Schur multipliers and pointing out the reference [Ka].

\section*{Notation}
We denote by $\bar n:= \{1, \dots, n\}$. Furthermore, we fix
$k$ to be a field of characteristic 0.
The reader can think of $\bf C$ if he or she wishes.
The theory is the same if $k$ is a super-commutative $\bf Q$ algebra and
(super-)vector spaces and dimensions are replaced by free modules and ranks.
Finally, if we fix a group $G$ then all remains true for a field of
a characteristic prime to $|G|$.

\section{Orbifold Frobenius algebras}
Recall the following definitions first presented in [K2] and
contained in [K3].

We fix a finite group $G$ and denote its unit element by $e$.

\subsection{Definition}
 A {\em G--twisted Frobenius algebra} (or $G$--Frobenius algebra for short) 
over
a field  $k$ of characteristic 0 is 
$<G,A,\circ,1,\eta,\varphi,\chi>$, where

\begin{tabular}{ll}
$G$&finite group\\
$A$&finite dim $G$-graded $k$--vector space \\
&$A=\oplus_{g \in G}A_{g}$\\
&$A_{e}$ is called the untwisted sector and \\
&the $A_{g}$ for $g \neq 
e$ are called the twisted sectors.\\
$\circ$&a multiplication on $A$ which respects the grading:\\
&$\circ:A_g \otimes A_h \rightarrow A_{gh}$\\
$1$&a fixed element in $A_{e}$--the unit\\
$\eta$&non-degenerate bilinear form\\
&which respects grading i.e. $\eta|_{A_{g}\otimes A_{h}}=0$ unless
$gh=e$.\\
\end{tabular}

\begin{tabular}{ll}
$\varphi$&an action by algebra automorphisms  of $G$ on $A$, \\
&$\varphi\in \mathrm{Hom}_{k-alg}(G,A)$, s.t.\
$\varphi_{g}(A_{h})\subset A_{ghg^{-1}}$\\
$\chi$&a character $\chi \in \mathrm {Hom}(G,k^{*})$ \\

\end{tabular}

\vskip 0.3cm         

\noindent satisfying the following axioms:

{\sc Notation.} We use a subscript on an element of $A$ to signify that it has homogeneous group 
degree  --e.g.\ $a_g$ means $a_g \in A_g$--, and we write $\varphi_{g}:= \varphi(g)$ and $\chi_{g}:= \chi(g)$. We also drop the subscript if $a\in A_e$.

\begin{itemize}

\item[a)] {\em Associativity}

$(a_{g}\circ a_{h}) \circ a_{k} =a_{g}\circ (a_{h} \circ a_{k})$
\item[b)] {\em Twisted commutativity}

$a_{g}\circ a_{h} = \varphi_{g}(a_{h})\circ a_{g}$
\item[c)]
{\em $G$ Invariant Unit}: 

$1 \circ a_{g} = a_{g}\circ 1 = a_g$

and

$\varphi_g(1)=1$
\item[d)]
{\em Invariance of the metric}: 

$\eta(a_{g},a_{h}\circ a_{k}) = \eta(a_{g}\circ a_{h},a_{k})$

\item[i)]
{\em Projective self--invariance of the twisted sectors}

$\varphi_{g}|A_{g}=\chi_{g}^{-1}id$

\item[ii)]
{\em $G$--Invariance of the multiplication}

$\varphi_{k}(a_{g}\circ a_{h}) = \varphi_{k}(a_{g})\circ  \varphi_{k}(a_{h})$

\item[iii)]

{\em Projective $G$--invariance of the metric}

$\varphi_{g}^{*}(\eta) = \chi_{g}^{-2}\eta$

\item[iv)]
{\em Projective trace axiom}

$\forall c \in A_{[g,h]}$ and $l_c$ left multiplication by $c$: 

$\chi_{h}\mathrm {Tr} (l_c  \varphi_{h}|_{A_{g}})=
\chi_{g^{-1}}\mathrm  {Tr}(  \varphi_{g^{-1}} l_c|_{A_{h}})$
\end{itemize}

An alternate choice of data is given by a one--form $\eps$, the co--unit 
with $\eps \in A_e^*$ and a three--tensor
$\la \cdot, \cdot, \cdot \ra \in A^* \otimes A^* \otimes A^*$  
which is of group degree $e$,
 i.e. $\la \cdot, \cdot, \cdot \ra|_{A_{g}\otimes A_{h}\otimes A_{k}}=0$ 
unless
$ghk=e$.

The relations between $\eta,\circ$ and $\eps,\mu$ 
are given by dualization. 

We denote by $\rho\in A_e$ the element dual to $\eps\in A_e^*$ and Poincar\'e dual to $1\in A_e$.

In the graded case, we call the degree $d$ of $\rho$ the degree of $A$.
This means that $\eta$ is homogeneous of degree $d$.

\subsection{Super-grading}
\label{super} We can enlarge the framework by considering
super--algebras rather than algebras. This will introduce the
standard signs.

The action of $G$ as well as the untwisted sector should be even.
The axioms that change are

\begin{itemize}

\item[b$^{\sigma}$)] {Twisted super--commutativity}

$a_{g}\circ a_{h} = (-1)^{\tilde a_g\tilde a_h} \varphi_{g}(a_{h})\circ a_{g}$

\item[iv$^{\sigma}$)]
{Projective super--trace axiom}

$\forall c \in A_{[g,h]}$ and $l_c$ left multiplication by $c$:

$\chi_{h}\mathrm {STr} (l_c  \varphi_{h}|_{A_{g}})=
\chi_{g^{-1}}\mathrm  {STr}(  \varphi_{g^{-1}} l_c|_{A_{h}})$
\end{itemize}
where $\mathrm{STr}$ is the super--trace.

Here we denoted by $\tilde a$ the $\mathbb{Z}/2\mathbb{Z}$ degree of $a$.

\subsection{$G$--graded tensor product} Given two $G$--Frobenius algebras
$\la G,A,\circ,1,\eta,\varphi,\chi\ra$ and
$\la G,A',\circ',1',\eta',\varphi',\chi'\ra$ we
defined \cite{K1} their tensor product
as $G$--Frobenius algebras to be the $G$--Frobenius algebra

$\la G,\bigoplus_{g \in G}( A_g \otimes A'_g),
\circ\otimes \circ',1\otimes 1',\eta\otimes \eta',\varphi\otimes \varphi',
\chi\otimes \chi'\ra$.

We will use the short hand notation $A \hat \otimes A'$ for this product.

\subsection{Grading and Shifts}
\label{shifts}
\subsubsection{Notation}
We denote by $\rho_g \in A_g$ the element defining $\eta_g$ and by
$d_g := \deg(\rho_g)$ the degree of $A_g$ and $s_g :=\deg(1_g)$ will 
be called the degree shift. We also set
$$
s^+_g  :=\frac{1}{2}(s(g) + s(g^{-1})) \quad s^- :=\frac{1}{2}(s(g) -s(g^{-1}))
$$
the degree defect.

Notice that $d= d_e$ if $d$ denotes the degree of $A$ given by $\eta$.

By considering $\eta|_{A_g \otimes A_{g^{-1}}}$ we find:

\subsubsection{Lemma} [K3] 
$$
s^+_g= d-d_g
$$

Notice there is no restriction (except anti--symmetry) on $s^-$.

The shift $s^-$ is not fixed, however, there is a standard choice
provided there exists a canonical choice of linear representation of $G$.

\subsubsection{Definition}
 The {\em standard shift} for a G--Frobenius algebra with a choice of linear
representation $\rho: G \rightarrow GL_n(k)$
is given by
$$s_{g}^+:= d-d_{g}$$
and
\begin{multline*}
s_{g}^- := \frac{1}{2\pi i}\mathrm{tr} (\log(g))-\mathrm{tr}(\log(g^{-1})):=
\frac{1}{2\pi i}(\sum_i \l_i(g)-\sum_i \l_i(g^{-1}))\\
=\sum_{i: \l_i \neq 0} (\frac{1}{2\pi i}2\l_i(g)-1)
\end{multline*}
where the $\l_i(g)$
are the logarithms of the eigenvalues
of $\rho(g)$ using the branch with arguments in $[0,2\pi)$ i.e.\
cut along the positive real axis.

In total we obtain:

$$
s_{g}= \frac{1}{2}(s_g^+ + s_g^-)= \frac{1}{2}(d-d_g)
+ \sum_{i:\l_i \neq0} (\frac{1}{2\pi i}\l_i(g)-\frac{1}{2})
$$

\subsubsection{Remark} This grading having its origin in physics
specializes to the so--called age grading or the orbifold grading of \cite{CR}
in the respective situations.
 
\subsection{Special $G$ Frobenius algebras}
\label{special}

\subsubsection{Definition} 
We call a $G$-Frobenius algebra {\em special} if all $A_g$ are cyclic
$A_e$ modules via the multiplication $A_e \otimes A_g \rightarrow A_g$
and there exists a collection of cyclic generators $1_g$ of $A_g$ such that
$\varphi_g(1_h)= \varphi_{g,h}1_{ghg^{-1}}$ with $\varphi_{g,h}\in k^*$.

The last condition is automatic, if the Frobenius algebra $A_e$ 
only has $k^*$ as invertibles, as is the case for cohomology algebras of
connected compact manifolds 
and Milnor rings of quasi--homogeneous functions with an
isolated critical point at zero.

Fixing the generators $1_g$ we obtain maps $r_g:A_e \rightarrow A_g$ by setting
$r_g(a_e)= a_e1_g$. This yields a short exact sequence

\begin{equation}
0\rightarrow I_g \rightarrow A_e \stackrel{r_g}{\rightarrow} A_g \rightarrow 0
\end{equation}

It is furthermore useful to fix a section $i_g$ of $r_g$.

We denote the concatenation $\pi_g:= i_g \circ i_g$.

\subsubsection{Special super $G$--Frobenius algebra}
The super version of special $G$--Frobe\-nius algebras is straightforward. Notice that
since each $A_g$ is a cyclic $A_e$--algebra 
its parity is fixed to be $(-1)^{\tilde g}:=\widetilde {1_g} $ times that of  $A_e$.
I.e.\ $a_g = i_{g}(a_{g})1_g$ and thus  $\widetilde{a_g} = \widetilde{i_g(a_{g})} \widetilde {1_g}$ .
In particular if $A_e$ is purely even $A_g$ is purely of degree $\tilde g$.

\subsubsection{Frobenius algebra structure on the twisted sectors}

Recall that the $A_{g}$ are Frobenius algebras
by the multiplication

\begin{equation}
a_g \circ_g b_g = i_g(a_g) i_g(b_g) 1_g
\end{equation}
and metric
\begin{equation}
\eta_g(a_g,b_g):= \eta(i_g(a_g)1_g,i_{g}(b_g)1_{g^{-1}})
\end{equation}

\subsubsection{Definition} Given a Frobenius algebra $A_e$ and a collection of 
cyclic $A_e$--modules $A_g:g \in G$
{\it a graded cocycle}  is a map $\g: G\times G \rightarrow A_e$
which satisfies 
$$\g(g,h)\g(gh,k)\equiv \g(g,hk)\g(h,k) \; \mathrm{ mod }\;  I_{ghk}$$

Such a cocycle is called
{\em section independent} if
$$(I_g + I_h)\g_{g,h} \subset I_{gh}$$
 
Two such cocycles are considered to be the same 
if $\g_{g,h} \equiv \g'_{g,h} \; \mathrm{ mod }\; I_{gh}$ and isomorphic, 
if they
are related by the usual scaling for group cocycles.

Given non--degenerate parings $\eta_g$ on the $A_g$, 
a cocycle is said to be {\em compatible with the metric}, if 
$$
\check r_g(1_g) = \g(g,g^{-1})
$$
where $\check r$ is the dual in the sense of vector spaces with
non--degenerate metric.

\subsubsection{The multiplication} 

Fixing a cyclic generator $1_g \in A_g$, the multiplication defines a 
section independent
graded cocycle $\g$ compatible which is compatible with the metric.
The cocycle $\g$ is defined via 
$$
1_g 1_h = \g_{g,h}1_{g,h}
$$

The section independence follows from the fact that 

$$(I_g+I_h)\g_{g,h}1_{gh}= (I_g+I_h)1_g1_h=0$$

In general, the multiplication is thus given by
\begin{equation}
\label{specialmult}
a_g b_h = i_g(a_g)i_h(b_h)\g_{g,h}1_{gh}
\end{equation}
for any choice of sections $i_g$.

The compatibility with the metric follows from the following equation
which holds for all $a\in A_e$:

$$
\eta(\g_{g,g^{-1}},a) = \eta(a1_{g},1_{g^{-1}})= \eta(r_{g}(a),1_{g^{-1}})
= \eta_{g}(1_{g},r_g(a))= \eta(\check r_g(1_g),a)
$$

\subsubsection{The $G$--action on the twisted sectors}
Consider a non--abelian cocycle $\varphi$ 
which is defined as a map $G\times G \rightarrow k^*$ satisfying:
\begin{equation}
\label{nonabcocycle}
\varphi_{gh,k}= \varphi_{g,hkh^{-1}}\varphi_{h,k}
\end{equation}
and
\begin{equation}
\varphi_{e,g}=\varphi_{g,e}=1
\end{equation}
where we used the notation $\varphi_{g,h}=\varphi(g,h)$

The $G$--action defines such a cocycle via
\begin{equation}
\varphi_g(1_h)= \varphi_{g,h}1_{ghg^{-1}}
\end{equation}
and in general the $G$--action is reduced to the one on the
non--twisted sector via
\begin{equation}
\varphi_g(a_h)=\varphi(g)(i_h(a_h)) \varphi_{g,h}1_{ghg^{-1}}
\end{equation}
for any choice of sections $i_h$.

\subsubsection{The compatibility equations}
The cocycles furthermore satisfy the following two compatibility equations:

\begin{equation}
\label{grpcompat}
\varphi_{g,h}\g_{ghg^{-1},g} = \g_{g,h} 
\end{equation}
and
\begin{equation}
\label{algaut}
\varphi_{k,g} \varphi_{k,h} \g_{kgk^{-1},khk^{-1}} 
= \varphi_{k} (\g_{g,h}) \varphi_{k,gh}
\end{equation}

We call a pair of a section independent cocycle and 
a non--abelian cocycle {\em compatible} if they satisfy the equations
(\ref{grpcompat}) and (\ref{algaut}).

\subsection{Definition}
A special $G$ reconstruction datum is a collection of Frobenius
algebras $(A_g,\eta_g,1_g): g\in G$ together with
an action of $G$ by algebra automorphisms
on $A_e$ and the structure of a cyclic $A_e$ module algebra on each $A_g$ with
generator $1_g$ such that $A_g$ and
$A_g^{-1}$ are isomorphic as  $A_e$ modules algebras.

\subsection{Theorem}(Reconstruction \cite{K2}) 
Given a special $G$ reconstruction datum
the structures of special $G$--Frobenius algebras are in 1--1
correspondence with compatible pairs of
a graded, section independent $G$ 2--cocycle with values in $A_e$ that is
compatible with the metric and a
non--abelian $G$ 2--cocycle with values in $K^{*}$,
satisfying the following conditions:
\begin{itemize}
\item[i)]$\varphi_{g,g}=\chi_g^{-1}$
\item[ii)]
$\eta_{e}(\varphi_{g}(a),\varphi_{g}(b)) =
    \chi_{g}^{-2}\eta_{e}(a,b)$
\item[iii)] The projective trace axiom
$\forall c \in A_{[g,h]}$ and $l_c$ left multiplication by $c$:
\begin{equation}
\chi_{h}\mathrm {Tr} (l_c  \varphi_{h}|_{A_{g}})=
\chi_{g^{-1}}\mathrm  {Tr}(  \varphi_{g^{-1}} l_c|_{A_{h}})
\end{equation}
\end{itemize}

\subsection{Rescaling}
Given a special G--Frobenius algebra, we can rescale the cyclic generators by
$\l_g$, i.e.\ we take the same underlying G--Frobenius algebra, but rescale the maps $r_g$ to $\tilde r_g$ with $\tilde 1_g=\tilde r_g(1)= \l_g 1_g$.
We also fix $\l_e=1$ to preserve the identity.

This yields an action of $\mathrm{Map}_{\text{pointed spaces}}(G,k^*)$
on the cocycles $\g$ and $\varphi$ preserving the underlying G--Frobenius algebra structure.

The action is given by:
\begin{eqnarray}
\label{rescaling}
\g_{g,h} &\mapsto& \tilde \g_{g,h}= \frac{\l_g \l_h}{\l_{gh}}\g_{g,h}\nn\\
\varphi_{g,h}&\mapsto&\tilde \varphi_{g,h}=
\frac{\l_h}{\l_{ghg^{-1}}}\varphi_{g,h}
\end{eqnarray}

\subsubsection{Remark} We can introduce the groups associated with the 
classes under this scaling and see that the classes of $\g$ correspond to classes in 
$H^2(G,A)$. We can also identify the non--abelian cocycles $\varphi$ with
one--group cocycles with values in $k^*[G]$ where we treat $k^*[G]$ as an 
abelian group with diagonal multiplicative composition
\begin{equation}
(\sum_g \l_g g)\cdot(\sum_h \mu_h h) := \sum_g \l_g\mu_g g
\end{equation}

and G--action given by conjugation:
\begin{equation}
s(g)(\sum_h \l_h h) = \sum_h \l_h ghg^{-1}
\end{equation}

This is done as follows:

We view the collection $\varphi_{g,.}$ as an element of $k^*[G]$
via

\begin{equation}
\varphi_g:=\sum_h \varphi_{g,h} ghg^{-1} 
\end{equation}

then 
$$
\varphi_{gh} =s(g) \varphi_h \cdot \varphi_g
$$

Indeed 
\begin{eqnarray*}
s(g) \varphi_h \cdot\varphi_g&=&
s(g)(\sum_k \varphi_{h,k}\, hkh^{-1})\cdot \sum_k \varphi_{g,k} \, gkg^{-1}\\ 
&=& \sum_k \varphi_{h,k} \, ghkh^{-1}g^{-1}\cdot \sum_k \varphi_{g,k} \, gkg^{-1}\\
&=& \sum_k \varphi_{h,k}\varphi_{g,hkh^{-1}}\,  ghkh^{-1}g^{-1}\\
&=& \sum_k \varphi_{gh,k}\, (gh)k(gh)^{-1}
\end{eqnarray*}

In this identification, 
equivalence under scaling corresponds to taking cohomology classes.

The trivial cocycles are of the form $s(g)a\cdot a^{-1}$ 
with $a = \sum \mu_g \, g$
\begin{equation}
s(g)a\cdot a^{-1}= \sum_h \mu_h ghg^{-1} \cdot \sum_h \mu_h^{-1} h
= \sum_h \frac{\mu_h}{\mu_{ghg^{-1}}} h
\end{equation}
and 
\begin{equation}
\tilde \varphi_g= \sum_h \tilde \varphi_{g,h} ghg^{-1}=
\sum_h \varphi_{g,h} ghg^{-1} \sum \frac{\l_h}{\l_{ghg^{-1}}} ghg^{-1}
= \varphi_g \cdot (s(g)a\cdot a^{-1})
\end{equation}
with $a = \sum_h \l_h h$.

It is clear that we could also take logarithms of the $\varphi$ 
and then we would
get cocycles with values in $k[G]$, but there is the problem 
of choosing a cut as it manifests itself 
in the setting of special G--Frobenius algebras in the definition 
of the degree shifts.

\subsubsection{Lemma}
\label{dimension}
Let  $A$ and $A_g$ be a graded Frobenius algebras with the top
degree of $A_g$ being $d_g$ then for a section independent cocycle
$\g_{g,g^{-1}} \subset L\subset A_e $ with $\dim(L)=\dim(A_g^{d_g})$,
where the superscript denotes a fixed degree.

{\bf Proof.}

By section independence
$$
I_{g}\g_{g,g^{-1}}=0
$$
Thus
$$
\g_{g,g^{-1}}\in  (i_g(A_g)^*)^{d-s_g^+}
$$
where ${}^*$ is the dual w.r.t.\ the form $\eta$ and we use the
splitting induced by the sections $i$ (N.B.\ if $\eta$ is also
positive definite, we could use an orthogonal splitting)
\begin{equation}
\label{splitting} A^{k}= I_g^{k}\oplus (i_g(A_g))^{k}
\end{equation}
and superscripts denote fixed degree. Furthermore

\begin{center}
$ \dim((i_g(A_g)^*)^{d_g}) = \dim(i_g(A_g)^{d_g}) = \dim(A^{d_g})-
\dim(I_{g})= \dim(A^{d_g})- \dim(\mathrm{Ker}(r_{g})|_{A^{d_g}})=
\dim(\mathrm{Im}(r_{g})|_{A^{d_g}})= \dim(A_{g}^{d_g})$
\end{center}
where we used the non--shifted grading on $A_{g}$. Thus
$\g_{g,g^{-1}}$ is fixed up to a constant.

If $\dim A_g=1$ then $\g_{g,g^{-1}}$ is fixed up to normalization
by the condition of section independence. The freedom to scale
$\g_{g,g^{-1}}$ is the same freedom one has in general for
choosing a metric for an irreducible Frobenius algebra. Recall
that in this case the space of invariant metrics is one
dimensional.

\subsection{Lemma} 
\label{dual}
If $a= i_g(a_g)\in i_{g}(A_{g})$ then 
$a\g_{g,g^{-1}} = \check r_{g}(a_g)$ and furthermore 
$i_{g}(A_{g})^*=\g_{g,g^{-1}} i_{g}(A_{g})$ where ${}^*$ is the 
Poincar\'e dual w.r.t.\ $\eta$ and the splitting (\ref{splitting}).
Moreover if $aI_g= 0$ then $a = \tilde a\g_{g,g^{-1}}$ for some
$\tilde a\in i_g(A_g)$.

{\bf Proof.}
For the first statement notice that:
$$
\eta(i_{g}(a_{g})\g_{g,g^{-1}},b)= \eta_{g}(a_{g},r_{g}(b))
$$
the second and third statement follow from this using the 
non--degenerate nature of $\eta, \eta_g$ and the splitting 
(\ref{splitting}). N.B.\ The statement is actually independent 
of the choice of splitting.

\subsection{Proposition}
\label{zerocheck}
If $\g_{g,h}=0$ then $\pi_{h}(\g_{g,g^{-1}})=0$ and
$\pi_{g}(\g_{h,h^{-1}})=0$

{\bf Proof.}

If $\g_{g,h}=0$ then\\
$0=\pi_{h}(\g_{g^{-1},gh}\g_{g,h})=
\pi_{h}(\g_{{g^{-1}},g}\g_{e,h}) = \pi_{h}(\g_{g^{-1},g}) = 
\pi_{h}(\g_{g,g^{-1}})$ 
and also

$0=\pi_{g}(\g_{g,h}\g_{gh,h^{-1}})=
\pi_{g}(\g_{g,e}\g_{h^{-1},h}) = \pi_{g}(\g_{h,h^{-1}}) $

\subsection{Definition}
\label{tranverse1} We call $A_g$ and $A_h$ transversal if
$s_{g}+s_{h}=s_{gh}$ and $s_{g^{-1}}+s_{h^{-1}}=s_{(gh)^{-1}}$.

From the section independence, we obtain:
\subsubsection{Lemma} If $A$ is irreducible and $A_g$ and $A_h$ 
are transversal and $\g_{g,h}\neq 0$ then 

$$I_g +I_h = I_{gh}$$

\subsection{Proposition} The converse of \ref{zerocheck} it true if
$A_g$ and $A_h$ are transversal.

{\bf Proof.} If $A_g$ and $A_h$ are transversal then
 $\deg(\g_{g,h})=0$ and $\g_{g,h}\in k$. 
The same holds for $\g_{h^{-1},g^{-1}}$. 
By associativity:
$$1_{g}1_{h}1_{h^{-1}}1_{g^{-1}}= \g_{h,h^{-1}} \g_{g,g^{-1}}
= \g_{g,h}\g_{h^{-1},g^{-1}}\g_{(gh),(gh)^{-1}}$$
and since $\g_{(gh),(gh)^{-1}}\neq 0$, we see that if $\g_{g,h}\neq 0$ and 
$\g_{h^{-1},g^{-1}}\neq 0$ then\\
 $ \g_{h,h^{-1}} \g_{g,g^{-1}}\neq 0$ 
so  $\pi_{h}(\g_{g,g^{-1}})\neq 0$ and $\pi_{g}(\g_{h,h^{-1}}) \neq 0$.

\subsection{Lemma} If $[g,h]=e$ 

\begin{equation}
\label{doubleconj}
\varphi_{g,h}= \varphi_{kgk^{-1},khk^{-1}}
\end{equation}

{\bf Proof.} 
$\varphi_{kgk^{-1},khk^{-1}}=\varphi_{k,h}
\varphi_{g,h}\varphi_{k^{-1},khk^{-1}}=
\varphi_{k,h}
\varphi_{g,h}\varphi^{-1}_{k,h}= \varphi_{g,h}$

\section{Discrete Torsion}\label{disc}
\subsection{The twisted group ring $k^{\a}[G]$}

Recall that given an element $\a \in Z^2(G,k^*)$
one defines the twisted group ring
$k^{\a}[G]$ to be given by the same linear structure with multiplication
given by the linear extension of

\begin{equation}
g\otimes h \mapsto \a(g,h) gh
\end{equation}
with $1$ remaining the unit element.
To avoid confusion we will denote elements of $k^{\a}[G]$ by
$\hat g$ and the multiplication with $\cdot$
Thus
$$
\hat g \cdot \hat h = \a(g,h) \widehat{gh}
$$

For $\a$ the following equations hold:
\begin{equation}
\a (g,e) = \a(e,g)=1, \qquad
\a(g,g^{-1})=\a(g^{-1},g)
\end{equation}
Furthermore
$$
\hat{g}^{-1}= \frac{1}{\a(g,g^{-1})}\widehat{g^{-1}}
$$
and
$$
\hat g\cdot \hat h\cdot\hat {g}^{-1} =
\frac{\a(g,h)\a(gh,g^{-1})}{\a({g,g^{-1})}}\widehat{ghg^{-1}}
= \frac{\a(g,h)}{\a(ghg^{-1},g)} \widehat{ghg^{-1}}=
\eps(g,h)\widehat{ghg^{-1}}
$$
with
\begin{equation}
\eps(g,h):=\frac{\a(g,h)}{\a(ghg^{-1},g)}
\end{equation}
\subsubsection{Remark}
\label{conorm} If the field $k$ is algebraically closed
we can find a representative for each class $[\a]\in H^2(G,k^*)$
which also satisfies

$$\a(g,g^{-1})=1$$

\subsubsection{Supergraded twisted group rings}
\label{supergroupring}
Fix $\a \in Z^2(G,k^*), \s \in \Hom(G,\Zz)$ then there is a twisted
super--version of the group ring where now the relations
read
\begin{equation}
 \hat g \hat h =  \a(g,h)\widehat {gh}
\end{equation}
and the twisted commutativity is
\begin{equation}
 \hat g \hat h = (-1)^{\s(g)\s(h)}\varphi_{g}(\hat h) \hat g
\end{equation}
and thus
\begin{equation}
 \varphi_{g}(\hat h)=
(-1)^{\s(g)\s(h)}\a(g,h)\a(gh,g^{-1}) \widehat{ghg^{-1}} =:
\varphi_{g,h} \widehat{ghg^{-1}}
\end{equation}
and thus
\begin{equation}
\eps(g,h) := \varphi_{g,h} = (-1)^{\s(g)\s(h)}\frac{\a(g,h)}{\a(ghg^{-1},g)}
\end{equation}

 We would just like to remark that the
axiom iv$^{\sigma})$ of  \ref{super} shows the difference between
super twists and discrete torsion.

\subsection{Definition}
We denote the $\a$-twisted group ring
with super--structure $\s$ by $k^{\a,\s}[G]$.
We still denote $k^{\a,0}[G]$ by $k^{\a}[G]$
where $0$ is the zero map and we denote $k^{0,\s}[G]$
just by $k^{\s}[G]$ where $0$ is the unit of the group $H^2(G,k^*)$.

A straightforward calculation shows
\subsection{Lemma}
 $k^{\a,\s}[G] = k^{\a}[G]\otimes k^{\s}[G]$.

\subsubsection{The $G$--Frobenius Algebra structure of $k^{\a}[G]$}
Fix  $\a \in Z^2(G,k^*)$.
Recall from \cite{K1,K2} the following structures which turn
$k^{\a}[G]$ into a special $G$--Frobenius algebra:

\begin{eqnarray}
 \g_{g,h}=\a(g,h) &&\eta(\hat g,\widehat{g^{-1}}) =\a(g,g^{-1})\nn\\
\chi_g= (-1)^{\tilde g}
&&\varphi_{g,h}=
\frac{\a(g,h)}{\a(ghg^{-1},g)}=:\eps(g,h)
\end{eqnarray}
\subsubsection{Relations}
\label{schurdt}
The $\eps(g,h)$ which are by definition given as
$\eps(g,h):= \frac{\a(g,h)}{\a(ghg^{-1},h})$ satisfy the equations:
\begin{eqnarray}
\eps(g,e)&=&\eps(g,g)=1\\\nn
 \eps(g_1g_2,h)&=&
\eps(g_1,g_2hg_2^{-1})\eps(g_2,h)\nn\\
 \eps(k,gh)& =& \eps(k,g)\eps(k,h)\frac{\a(kgk^{-1},khk^{-1})}{\a(g,h)}\nn\\
\eps(h,g)&=&\eps(g^{-1},ghg^{-1})
\frac{\a([g,h],h)}{\a([g,h],hgh^{-1})}
\end{eqnarray}

This yields for {commuting elements}:

\begin{eqnarray}
\label{eps}
\eps(g,e)=\eps(g,g)=1 && \eps(g,h) =\eps(h^{-1},g)=\eps(h,g)^{-1}\nn\\
\eps(g_1g_2,h)= \eps(g_1,h)\eps(g_2,h)  &&
\eps(h,g_1g_2) = \eps(h,g_1) \eps(h,g_2)
\end{eqnarray}

In the physics literature discrete torsion is sometimes defined to
be a function  $\eps$ defined on commuting elements of $G$ taking
values in $U(1)$ and satisfying the equations (\ref{eps}).

\subsection{The trace axiom}
The trace condition for
 non--commuting elements  reads
$$
 (-1)^{\tilde h} (-1)^{\tilde g}\varphi_{h,g}
\g_{[g,h],hgh^{-1}} 
=
 (-1)^{\tilde g} (-1)^{\tilde h}\varphi_{g^{-1},ghg^{-1}}
\g_{[g,h],h}
$$
stripping off the sign, we  rewrite the l.h.s.\ as
\begin{eqnarray*}
\varphi_{h,g}
\g_{[g,h],hgh^{-1}} \widehat{gh}
&=&\varphi_{h,g}\g_{g,h}^{-1}\widehat{[g,h]}\widehat{hgh^{-1}} \widehat{h} \\
= \varphi_{h,g}\g_{[g,h],hg}\g_{hgh^{-1},h} \g_{g,h}^{-1}\widehat{gh}
&=&\g_{[g,h],hg}
\g_{h,g}
\g_{g,h}^{-1}\widehat{gh}
\end{eqnarray*}
and the r.h.s.\ can be rewritten as
\begin{eqnarray*}
\varphi_{g^{-1},ghg^{-1}}
\g_{[g,h],h} \widehat{gh}&=&
\varphi_{g^{-1},ghg^{-1}} \g_{ghg^{-1},g}^{-1}
\widehat{[g,h]}\widehat{h} \widehat{g}\\ 
=\varphi_{g^{-1},ghg^{-1}}\g_{[g,h],hg}\g_{h,g} \g_{ghg^{-1},g}^{-1}
\widehat{gh}
&=&\g_{[g,h],hg}\g_{h,g}\g_{g,h}^{-1} \widehat{gh}
\end{eqnarray*}
which coincides with the calculation above.

This is of course all clear if $[g,h]=e$, but there is no restriction that
the group be commutative.

\subsubsection{Remark} The function $\eps$ can be interpreted as
a cocycle in $Z^1(G,k^*[G])$ where $k^*[G]$ are the elements of
$k[G]$ with invertible coefficients regarded as a $G$ module by
conjugation (cf. \cite{K1,K2}). This means in particular that on
{\em commuting elements} $\eps$ only depends on the class of the cocycle $\a$.

\subsection{Theorem} The possible super $G$ Frobenius algebra
structures on $A=\bigoplus_{g\in G} k$
are the structures of super twisted group rings.
The isomorphism classes of these algebras correspond to pairs of a class
$[\a]\in H^2(G,K^*)$ and a homomorphism $\sigma \in {\rm Hom}(G,\Zz)$.

{\bf Proof.} Assume that we have a $G$ Frobenius algebra structure on $A$ then it
is a special $G$--Frobenius algebra since $1\in A_e$ is the unit.
Then due to the non--degeneracy of the metric $\g_{g,g^{-1}}\in k^*$ furthermore
$\pi_h(\g_{g,g^{-1}})=\g_{g,g^{-1}})\in k^*$ and thus by \ref{zerocheck} 
$\forall g,h \in G:\g_{g,h}\in k^*$, thus $\g \in Z^2(G,k^*)$ and by compatibility
the $\varphi$ are fixed. Lastly, since $ \g_{g,h}\in k^*$ and $\tilde\g_{g,h} =0$
the supergrading $\tilde{}$ must be a homomorphism, i.e.\ $\tilde{}\in {\rm Hom}(G,\Zz)$.

Vice versa the construction above shows that given a cycle $\a \in Z^2(G,k^*)$
and a homomorphism $\sigma \in {\rm Hom}(G,\Zz)$ we get a structure of 
super $G$ Frobenius algebra with the underlying data.

The statement about the isomorphisms classes follows directly from rescaling.

\subsection{The action of discrete Torsion}

\subsubsection{The action of $Z^2(G,k^*)$ }
The group $Z^2(G,k^*)$ acts naturally on 
$Z^2(G,A)$ via $(\a,\g)\mapsto \g^{\a}:=\g\cdot\a$ and
on $H^1(G,k^*[G])$ via $(\a,\varphi)\mapsto \varphi^{\a}:= \eps_{\a}\cdot \varphi$
were $\eps_{\a}(g,h)=\frac{\a(g,h)}{\a(ghg^{-1},g})$.

We call this action action by $\a$ twist or by the discrete torsion 
$\a$.

\subsection{Definition}
Given a $G$--Frobenius algebra
$A$ and an element $\a \in Z^2(G,k)$, we define the
$\a$--twist (or the twist by the discrete torsion $\a$) 
of $A$ to be the $G$--Frobenius algebra
$A^{\a}:= A \hat\otimes k^{\a}[G]$.

\subsection{Proposition}
\label{defprop}
Notice that as vector spaces
\begin{equation}
\label{alphaiso}
A^{\a}_{g}= A_g \otimes k \simeq A_g
\end{equation}
Using this identification the $G$--Frobenius structures  given by
(\ref{alphaiso}) are
\begin{eqnarray}
\circ^{\a}|_{A^{\a}_{g}\otimes A^{\a}_{h}}= \a(g,h) \circ &&
\varphi^{\a}_g|_{A^{\a}_h}=\eps(g,h)\varphi_g\nn\\
\eta^{\a}|_{A^{\a}_g\otimes A^{\a}_{g^{-1}}}= \a(g,g^{-1})\eta&&
\chi_g=\chi_g
\end{eqnarray}

\subsection{Lemma}
Let $\la G,A,\circ,1,\eta,\varphi,\chi\ra$  be
a $G$--Frobenius algebra or more generally a super Frobenius algebra with
super grading $\tilde{}\;\in {\rm Hom}(A,\Zz)$
 then
$A\otimes k^{\s}[G]$ is isomorphic to the super $G$--Frobenius algebra
$\la G,A,\circ^{\sigma},1,\eta^{\sigma},\varphi^{\s},\chi^{\s}\ra$ with 
super grading
${}^{\sim \s}$, where
\begin{eqnarray*}
\circ^{\sigma}|_{A_g\otimes A_h}=(-1)^{\tilde g \sigma(h)}\circ
&\quad &
\varphi^{\s}_{g,h} = (-1)^{\s(g)\s(h)}\varphi_{g,h}\\
\eta_g^{\sigma}=(-1)^{\tilde g \sigma(g)}\eta_g
&\quad& \chi^{\s} = (-1)^{\s(g)}\chi_g\\
\tilde a_g^{\s}= \tilde a_g + \s(g)&&\\
\end{eqnarray*}

\subsection{Definition}
Given a $G$--Frobenius algebra
$A$ a twist for $A$ is a pair of functions
$(\l:G\times G \rightarrow k^*,\mu:G\times G \rightarrow k^*)$
such that $A$ together with the new $G$--action
$$\varphi^{\l}(g)(a) = \oplus_h \l(g,h) \varphi(g)(a_h)$$
and the new multiplication
$$
a_g \circ^{\mu} b_h = \mu(g,h) a_g \circ b_h
$$
is again a $G$--Frobenius algebra.

A twist is called universal if it is defined for all $G$--Frobenius algebras.

\subsubsection{Remark} We could have started from a pair of functions
$(\l:A\times A \rightarrow k^*,\mu:G\times A \rightarrow k^*)$ in order
to projectively change the multiplication and $G$ action, but it is
clear that the universal twists (i.e.\ defined for any $G$--Frobenius
algebra) can only take into account
the $G$ degree of the elements.

\subsubsection{Remark}
These twists arise from a projectivization of the $G$--structures
induced on a module over $A$ as for instance the associated
Ramond--space (cf.\ \cite{K1}). In physics terms this means that
each twisted sector will have a projective vacuum, so that fixing
their lifts in different ways induces the twist. Mathematically
this means that the $g$ twisted sector is considered to be a Verma
module over $A_g$ based on this vacuum.

\subsection{Theorem}\cite{K4}
Given a (super) $G$--Frobenius algebra $A$
the universal twists are in 1--1 correspondence with
elements $\a \in Z^2(G,k^*)$  and the isomorphism classes of universal
twists are given by $H^2(G,k^*)$. Furthermore
the universal super  re--gradings are
in 1-1 correspondence with $\Hom(G,\Zz)$ and these
structures can be realized by tensoring with $k^{\s}[G]$
for $\s \in  \Hom(G,\Zz)$.

Here a super re--grading is a new super grading on $A$ with which
$A$ is a super $G$--Frobenius algebra and universal means that
the operation of re--grading is defined for all $G$--Frobenius algebras.

We call the operation of forming a tensor product with $k^{\a}[G]:
\a \in Z^2(G,k^*)$ a twist by discrete torsion. The term discrete
refers to the isomorphism classes of twisted $G$--Frobenius
algebras which correspond to classes in $H^2(G,k^*)$. Furthermore,
we call  the operation of forming a tensor product with
$k^{\s}[G]:\s \in {\rm Hom}(G,\Zz)$ super--twist.

\subsection{Remark} If $k$ is algebraically closed, then in each class
of $H^2(G,k^*)$ there is a representative with $\a(g,g^{-1})=1$.
Using these representatives it is possible to twist a special
$G$--Frobenius algebra without changing its underlying special
reconstruction data.

\section{Functorial setup}

The functorial setup of orbifold Frobenius algebras and reconstruction 
is discussed in the following.

Let $\Frob$ be the category of Frobenius algebras, whose objects are Frobenius 
algebras and morphisms are morphisms which respect all the structures.

\subsection{Definitions} 

{\em A  $G$--category}
is a category $\C$ where for each object
$X \in Ob(\C)$ and each $g \in G$ there exists an object $X^g$ and a
 morphism
$i_g \in \Hom(X^g,X)$  with $X^e=X$ and $i_e= id$
and there are isomorphisms $\psi_{g,g^{-1}}\in \Hom(X^{g},X^{g^{-1}})$.

We call a category {\em a $G$ intersection category}
if it is a $G$ category and for each pair
$(g,h) \in G\times G$ and object $X\in Ob(\C)$ there are isomorphisms
$\psi\in \Hom( (X^g)^h, (X^h)^g)$ and morphisms
$i^{gh}_{g,h} \in \Hom((X^g)^h,X^{gh})$.

A {\em $G$--action} for a $G$--category is given by
a collection of morphisms $\phi_g(X,h) \in \Hom(X^h,X^{ghg^{-1}})$
which are compatible with the structural morphisms and
satisfy $\phi_g(X,g'hg^{\prime -1})\phi_{g'}(X,h)= \phi_{gg'}(X,h)$.

\subsection{Examples}
Examples of an intersection $G$--category with $G$--action are
categories of spaces equipped  with a $G$--action whose fixed
point sets are in the same category. Actually this is the category
of pairs $(X,Y)$ with $X$ say a smooth space with a $G$--action
and $Y$ a subspace of $X$. Then $(X,Y)^g:=(X,Y\cap Fix(g,X))$ with
$Fix(g,X)$ denoting the fixed points of $g\in G$ in X, and
$i_g=(id,\iota_g)$ with $\iota_g: Y\cap Fix(g,X)\rightarrow Y)$
being the inclusion. It is enough to consider pairs $(X,Y)$ where
$Y\subset X$ is the set fixed by a subgroup generated by an
arbitrary number of elements of $G$: $H:=\langle
g_1,\dots,g_k\rangle$

We could also consider the action on the $X^g$ to be trivial and
set $(X^g)^h:= X^g$. This will yield a $G$--category.

Also the category of functions $f:{\bf C}^n\rightarrow {\bf C}$
with an isolated singularity at $0$ together with a group action
of $G$ on the variables induced by a linear action of $G$ on the
linear space fixing the function is an example of a $G$--category.
This is a category of triples $({\bf C}^n, f: {\bf C}^n
\rightarrow {\bf C}), \rho \in {\rm Hom}(G,GL(n))$ such that $f$
has an isolated singularity at zero and $f(\rho({\bf z}))= f({\bf
z})$ for ${\bf z}\in {\bf C}^n$ with morphisms being linear
between the linear spaces such that all structures are compatible.
The functor under consideration is the local ring or Milnor ring.
Again we set $(X^g)^h:= X^g$.

 Here the role of
the fixed point set is played by the linear fixed point set and
the restriction of the function to this fixed point set
(cf.\cite{K1}). Again we can consider pairs of an object and a
subobject as above in order to get an intersection $G$--category.

Our main examples are smaller categories such as a global
orbifold. As a $G$ category, the objects are the fixed point sets
of the various cyclic groups generated by the element of $G$ and
the morphisms being the inclusion maps. Again we set $(X^g)^h:=
X^g$. For a global orbifold, we can also consider all fixed point
sets of the groups generated by any number of elements of $G$ as
objects together with the inclusion maps as morphisms. This latter
will render a $G$--intersection category.

The same is true for isolated singularities. Here the objects are
the restriction of the function to the various subspaces fixed by
the elements of $g$ together with the inclusion maps or for the $G$--intersection
category we consider all intersections of these subspaces together with the
restriction of the function to these subspaces as objects,
again with the inclusion morphisms.

Now, suppose we have a
$G$--category $\C$ and a contravariant functor $\mathcal{F}$
from $\C$ to $\Frob$. In this setting
there might be several schemes to define a ``stingy geometry'' by
augmenting the functor to take values in $G$--Frobenius algebras.
But all of these schemes have to have the same additive structure provided
by the ``classical orbifold picture'' (see \ref{classical}) and satisfy
the axioms of $G$--Frobenius algebras (see \S 2). Furthermore there
are more structures which are already fixed in this situation,
which is explained below. These data can sometimes be used to classify
the possible algebra structures and reconstruct it when the classification
data is known. In the case of so--called special $G$--Frobenius algebras
a classification in terms of group cohomology classes is possible.

There are some intermediate steps which contain partial information
that have been previously considered, like
the additive structure, dimensions etc., as discussed in \ref{classical}.

\subsubsection{The ``classical orbifold picture''}
\label{classical} Now, suppose we have a $G$--category $\C$ and a
contravariant functor $\mathcal{F}$ from $\C$ to $\Frob$, then for
each $X \in Ob(\C)$, we naturally obtain the following collection
of Frobenius algebras: $(\F(X^g):g\in G)$ together with
restriction maps $r_g = \F(i_g): \F(X) \mapsto \F(X^g)$.

One possibility is to regard the direct sum of the Frobenius algebras
$A_g:=\F(X^g)$.

The first obstacle is presented in the presence of a
grading, say by ${\bf N}, {\bf Z}$ or ${\bf Q}$; as it is well known that
the direct sum of two graded Frobenius algebras is only well defined
if their Euler dimensions  (cf.\ e.g.\ \cite{K3}) agree. This can, however,
be fixed by using the shifts
$s^+$ discussed in \ref{shifts}. If the grading was originally in ${\bf N}$
these shifts are usually in $\frac{1}{2}{\bf N}$, but in the complex case still lie in ${\bf N}$.

Furthermore, if we have a $G$--action on the $G$ category, it
will induce the structure of a $G$--module on this direct
sum.

Each of the Frobenius algebras $A_g$ comes equipped with
its own multiplication,
so there is a ``diagonal'' multiplication
for the direct sum which is the direct sum of these multiplications.

Using the shift $s^+$ it is possible to define a ``classical
theory'' by considering the diagonal algebra structure and taking
$G$--invariants. This is the approach used in  \cite{AS}, \cite{T}
and \cite{AR}. The paper \cite{AS} shows that this structure
describes the $G$--equivariant rather than the $G$--invariant
geometry.

One can of course forget
the algebra structure altogether and retain only
the additive structure. This was done e.g.\ in \cite{S} in the setting of
V--manifolds (i.e.\ orbifolds).
Concentrating only on the dimensions one arrives for instance at
the notion of ``stringy numbers'' \cite{BB}.

\subsubsection{The ``stringy orbifold picture''}

The ``diagonal'' multiplication is however {\em not} the right object to study
from the perspective of ``stringy geometry'' or a TFT with a finite
gauge group \cite{K1,CR}.
The multiplication should rather be $G$--graded, i.e.\ map
$A_g \otimes A_h \rightarrow A_{gh}$. We call such a product ``stringy'' product.

Here the natural question is the following:

{\bf Question.}
Given the additive structure
of a $G$--Frobenius algebra, what are the possible ``stringy'' products?

A more precise version of this question is the setting of our reconstruction program \cite{K2, K3}.

\subsubsection{The $G$--action} One part of the structure of 
a $G$--Frobenius algebra is the $G$--action. If the $G$--category
is already endowed with a $G$--action we can use it to reconstruct the
$G$--action on the $G$--Frobenius algebra, which in turn limits the
choices of ``stringy'' products to those that are compatible.

\subsubsection{Invariants} By definition $G$--Frobenius algebras
come with a $G$ action whose invariants form a commutative algebra.
Due to the nature of the $G$ action
this commutative algebra is graded by conjugacy classes, and under
certain conditions the metric descends and
the resulting algebra is again Frobenius. The induced multiplication
is multiplicative in the conjugacy classes and we call such
a multiplication commutative ``stringy''.

\subsubsection{Examples}
Examples of commutative ``stringy'' products
are orbifold (quantum) cohomology \cite{CR}.
For cohomology of global orbifolds
 it was shown in \cite{FG} and recently in \cite{JKK} that
there is a group graded version for global orbifold cohomology
which has the structure of a $G$ Frobenius algebra, as we had previously
postulated \cite{K2}.
For new developments on quantum
deformations of the $G$--Frobenius algebras see \cite{JKK}.

\subsubsection{Special $G$--Frobenius algebras}
The special reconstruction data reflects this situation in the special case
that the $A_g$ algebras are cyclic $A_e$ modules. This is a restriction 
which leads to an answer in terms of cocycles for a large
class of examples. This class includes all Jacobian Frobenius algebras as 
well as symmetric products and special cases of geometric actions on 
manifolds.

The general idea can be generalize to
non--cyclic case although computations get more involved.

\subsection{Definition}
Given a $G$--category $\C$, we call the tuple $(X^g):g\in G$ a G--collection.

The category of $G$--collections of a $G$--category is the category whose 
objects are $G$--collections and whose morphisms are collections of morphisms
$(f^g)$
s.t.\ the diagrams
$$\begin{matrix}
X^g&\stackrel{i_g}{\rightarrow}&X\\
\downarrow f^g&&\downarrow f\\
Y^g& \stackrel{i_g}{\rightarrow} & Y
\end{matrix}$$
commute.

\subsection{Definition} 
A G--Frobenius functor is a functor from the category of $G$--collections
of a $G$--category to $G$--Frobenius algebras.

\subsection{Reconstruction/classification}
The main question of the reconstruction/classification
 program is whether one can extend
a functor from a $G$--category $\C$ to Frobenius algebras to a
$G$--Frobenius functor, and if so how many ways are there to do this.

One can view this as the analogue of solving the associativity equations for
general Frobenius algebras. Some of the solutions correspond to
quantum cohomology, some to singularities, etc. and maybe others
to other ``string''--schemes. The structures of possible
``stringy'' products provide a common approach. The systematic
consideration of all possible products confines the choices of
string equivalents of classical concepts and allows to identify divers
approaches.

The answer to the main
question of reconstruction/classification can be
answered in the special case where all of the twisted
sectors are cyclic in terms of group cohomological data (see below).
This is the content of the Reconstruction Theorem of \cite{K1}.

The consequences are sometimes quite striking as in the case
of symmetric products, where there is only {one}
possible ``stringy'' orbifold
product.

The restrictions on the possible multiplicative structures are even
stricter if one is considering data stemming 
from a $G$--intersection category.

This is the content of the next section.

\section{Intersection G--Frobenius algebras}
We will now concentrate on the situation of functors from $G$--intersection
categories to Frobenius algebras.

Given a $G$--class in such a category a functor to Frobenius algebras will
provide the following structure which reflects the possibility
to take fixed point sets iteratively. Say we look at the fixed
points with respect to
elements $g_1, \dots, g_n$. These fixed point sets
will be invariant under the group spanned by the elements
$g_1, \dots, g_n$ and they are just the intersection of
the respective fixed point sets of the elements $g_i$.
The underlying spaces are therefore invariant with respect
to permutation of the elements $g_i$, and if $g$ appears twice
 among the $g_i$ then one can shorten
the list by omitting one of the $g_i$. Also if a list $g_i$
includes $g^{-1}$ we may replace it by $g$. Finally, the fixed
point set under the action of the group generated by two elements
$g$ and $h$ is a subset of the fixed point set of the group
generated by their product $gh$. Translating this into the
categorical framework, we obtain:

\subsection{Definition}A $G$--intersection Frobenius datum of level $k$ is the following:
For each collection $(g_1,\ldots, g_n)$ with $n\leq k$ of elements
of $G$, a Frobenius algebra $A_{g_1,\dots,g_n}$ and the following
maps:

Isomorphisms
$$\Psi_{\sigma}:
 A_{g_1,\dots,g_n}\rightarrow A_{g_{\sigma(1)},\dots,g_{\sigma(n)}}$$
for each $\sigma \in \Sn$ called {permutations}.

Isomorphisms
$$\Psi^{g_1,\dots, g_i, \dots, g_n}_{g_1,\dots, g_i^{-1},\dots ,g_n  }:
A_{g_1,\dots, g_i, \dots, g_n} \rightarrow
 A_{g_1,\dots, g_i^{-1},\dots ,g_n}$$
 commuting with the permutations.

Morphisms
$$r_{g_1,\dots, g_i, \dots g_n}^{g_1,\dots , \hat g_i,\dots ,g_n}:
A_{g_1,\dots,\hat g_i, \dots, g_{n}}
\rightarrow A_{ g_1,\dots,g_n}$$
commuting with the permutations. (Here the symbol $\hat {}$
is used to denote omission.)
Such that the diagrams
$$
\begin{matrix}
A_{g_1,\dots,\hat g_i, \dots, \hat g_j,\dots, g_{n}}
&\stackrel
{r_{g_1,\dots, g_i, \hat g_j,\dots g_n}^
{g_1,\dots , \hat g_i,\dots, \hat g_j, \dots ,g_n}}
{\rightarrow}& A_{ g_1,\dots,\hat g_j,\dots,g_n}\\
\downarrow 
r^{g_1,\dots,\hat g_i, \dots, \hat g_j,\dots, g_{n}}_{g_1,\dots, \hat g_j,
\dots, g_{n}}&&
\downarrow r^{ g_1,\dots,\hat g_j,\dots,g_n}_{g_1,\dots,g_n}\\
A_{g_1,\dots, \hat g_i,\dots, g_{n}}
&\stackrel{r^{g_1,\dots, \hat g_i,\dots, g_{n}}_{g_1,\dots,g_n}}
{\rightarrow}& A_{ g_1,\dots,g_n}
\end{matrix}
$$
are co--Cartesian.

Isomorphisms
$$i_{g_1,\dots, g, \dots , g,\dots, g_n}^{g_1,\dots g, \dots, \hat g,\dots ,g_n}:
A_{g_1,\dots, g, \dots , g,\dots, g_n}
\rightarrow A_{g_1,\dots g, \dots, \hat g,\dots ,g_n}$$
commuting with the permutations.

And finally morphisms:
$$r^{g_1,\dots,g_{i}g_{i+1},\dots, g_n}_{g_1,\dots,g_{i},g_{i+1},\dots, g_n}:
A_{g_1,\dots,g_{i}g_{i+1},\dots, g_n} \rightarrow
A_{g_1,\dots,g_{i},g_{i+1},\dots, g_n}$$
commuting with the permutations.

If this data exists for all $k$ we call the data simply
$G$--intersection Frobenius datum.

\subsection{Notation}
We set $r_{g_1,\dots,g_n}:= r_{g_1, \dots, g_n}^{g_1,\dots, g_{n-1}}
\circ \dots \circ r_{g_1}$
and we set $I_{g_1,\dots,g_n}:= \mathrm{Ker}( r_{g_1,\dots, g_n})$.
Notice that this definition of $I_{g_1,\dots,g_n}$ is independent
of the order of the $g_i$.

\subsection{Remarks}
\begin{itemize}
\item[1)]
In order to (re)--construct
a suitable multiplication on $\bigoplus A_g$ it is often
convenient to use the double and triple intersections (i.e.\  level 3).
Where the double intersection are used for the multiplication and
triple intersections are used to
show associativity.
\item[2)] We can use the double intersections to define $G$--Frobenius
algebras based on each of the $A_g$ i.e.\ on
 $\bigoplus_{h\in Z(g)} A_{g,h}$ for each fixed
$g$--where $Z(g)$ denotes the centralizer of $g$.
\end{itemize}

\subsubsection{Definition}
A {\em $G$--action} for an intersection $G$--Frobenius algebra of 
level $k$ is given by 
a collection of morphisms 
$$\phi_g(A_{g_1,\dots,g_n},h) \in 
\Hom(A_{g_1,\dots,g_n,h},A_{g_1,\dots,g_n,ghg^{-1}})$$
which are compatible with the structural homomorphisms and
satisfy 
$$\phi_g(A_{g_1,\dots,g_n},g'hg^{\prime -1})\phi_{\g'}
(A_{g_1,\dots,g_n},h)= \phi_{gg'}(A_{g_1,\dots,g_n},h)$$

\subsection{Definition} 
We call an intersection $G$ Frobenius datum
a special $G$ intersection Frobenius
datum, if all of the $A_{g_1,\dots,g_n}$ are cyclic $A_e$
module algebras via the restriction maps such that the $A_e$ module
structures are compatible with the restriction morphisms $r$. Here
the generators are given by $r_{g_1,\dots, g_n}(1)$ and the $A_e$
module structure is given by $a\cdot b:= r_{g_1,\dots,g_n}(a)b$.

\subsection{Remark}
In the case of special $G$--Frobenius algebras, the presence of
special intersection data gives a second way to look at
the multiplication. The first  way is to use the
restrictions $r_g$ and sections $i_g$ to define the multiplication
as discussed in \S \ref{special} (see eq.\ (\ref{specialmult})). A second possibility is to use
the intersection structure. This can be done in the following way:
first push forward to double intersections, second use the
Frobenius algebra structure there to multiply, then pull the
result back up to the invariants of the product, but allowing to
multiply with an obstruction class before pulling back. This is
discussed below in \S \ref{multiplication}.

The precise relation between the two procedures is given by
the following Proposition and \ref{specialmult}.
\subsection{Proposition}
\label{intersect}
 Given a special $G$ intersection datum
(of level $2$),
the following decomposition holds for section independent cocycles
$\g$:
\begin{equation}
r_{gh}(\g_{g,h}) =\check r_{g,h}^{gh}(\tilde\g_{g,h})
= i_{g,h}^{gh}(\tilde \g_{g,h})  \check r_{g,h}^{gh}(1_{g,h})
=\bar \g_{g,h}  \g_{g,h}^{\perp}
\end{equation}
 for some 
section $i_{g,h}$ of $r_{g,h}$, $\tilde \g_{g,h} \in (A_{g,h})^{e}$,
$\bar \g_{g,h}\in i_g,h)(A_g,h)$ of degree $e$.
and $ \g_{g,h}^{\perp}:=\check r_{g,h}^{gh}(1_{g,h})$.
Here $e=s_{g}+s_{h}-s_{gh}-s^{+}_{g,h}+s^{+}_{gh}$ with $s^{+}_{g,h} :=
d-d_{g,h}$ and $d_{g,h}=\deg(\rho_{g,h})$
and we again used the unshifted degrees. (In particular if  the
$s^{-}=0$ then $e= \frac{1}{2}(s_g^++s_h^++s_{gh}^+)-s_{g,h}^+
=\frac{1}{2}(d-d_g-d_h-d_{gh})+d_{g,h}$)

{\bf Proof.}
We notice that $I_{g}+I_{h} = I_{g,h}$ 
and $(I_{g}+I_{h})\g_{g,h}\subset I_{gh}$, and set $J:= r_{gh}(I_{g,h})$.
Choosing some section
$i_{g,h}^{gh}$ of $r_{g,h}^{gh}$, we can define the splitting
\begin{equation}
\label{split2}
A_{gh}^k = i_{g,h}^{gh}(A_{g,h})\oplus J
\end{equation}
where again ${}^{k}$ means the homogeneous component of degree $k$.
Now
$$
\g_{g,h}\in (i_{g,h}^{gh}(A_{gh})^*)^{e}
$$
where ${}^*$ is the dual w.r.t.\ the form $\eta_{gh}$ and the splitting
(\ref{split2}) and $e=s_{g}+s_{h}-s_{gh}+s^+_{gh}-s^+_{g,h}$.

From which the claim follows by an argument completely
analogous to the proof of Lemmas \ref{dimension} and
\ref{dual}.

Also generalizing the fact that 
\begin{equation}
I_g \g_{g} = I_g \check r_g (1_g)=0
\end{equation}
we obtain

\subsection{Lemma}
\begin{equation}
(I_g +I_h)\g_{g,h}^{\perp}\subset I_{g,h}
\end{equation}

\subsection{Multiplication} 
\label{multiplication}
From the section independence of $\g$, we
see for a special $G$--Frobenius algebra which is part of 
a special $G$--intersection 
Frobenius datum of level $\geq 2$ that the multiplication 
$A_g \otimes A_h \rightarrow A_{gh}$ can be factored through 
$A_{g,h}$. To be more precise, we have the following commutative diagram.

$$
\begin{matrix}
\label{multdiagram}
A_g\otimes A_h &\stackrel{\mu}{\rightarrow} &A_{gh}\\
\downarrow r^g_{g,h}\otimes r^h_{g,h}&&
\uparrow \check r^{g,h}_{gh}\circ l_{\tilde\g_{g,h}}\\
A_{g,h}\otimes A_{g,h}&\stackrel{\mu}{\rightarrow}&A_{g,h}
\end{matrix}
$$
where $l_{\tilde\g_{g,h}}$ is the left multiplication with $\tilde\g_{g,h}$.
That is  using the multiplication in $A_{g,h}$

\begin{equation}
a_g \circ b_h =
\check r^{gh}_{g,h}(r^g_{g,h}(a_g)r^h_{g,h}(b_h)\tilde\g_{g,h})
\end{equation}

\subsubsection{Remark} The decomposition into the terms
$\tilde \g$ and $\g^{\perp}$ can be understood as decomposing the
cocycle into a part which comes from the normal bundle of $X^{g,h}
\subset X^{gh}$ which is captured by $\g^{\perp}$ and an additional
obstruction part.

\subsection{Associativity equations}
\label{ass}
Furthermore in the presence of
a special $G$ intersection Frobenius datum of level $\geq 3$ the
 associativity equations can be factored through
$A_{g,h,k}$. More precisely, we have the following commutative diagram
of restriction maps:

\begin{equation}
\label{assdiagram}
\begin{matrix}
&&&&A_{ghk}&&&&\\
&&&\swarrow&&\searrow&&&\\
A_{gh}&\rightarrow&A_{gh,k}&&\downarrow&&A_{g,hk}&\leftarrow&A_{hk}\\
\downarrow&&&\searrow&&\swarrow&&&\downarrow\\
A_{g,h}&&\rightarrow&&A_{g,h,k}&&\leftarrow&&A_{h,k}\\
\end{matrix}
\end{equation}

More technically:
Using the associativity equations for the $\g$, we set
\begin{equation}
\label{tildetripel}
r_{ghk}(\g_{g,h} \g_{gh,k}):= 
\g_{g,h,k}
\end{equation}
and associativity says that also 
\begin{equation}
r_{ghk}(\g_{h,k} \g_{g,hk})= 
\g_{g,h,k}
\end{equation}

By analogous arguments as utilized above one finds 
\begin{equation}
\g_{g,h,k}= i_{g,h,k}^{ghk}(\tilde\g_{g,h,k}) \check r_{g,h,k}^{ghk}(1_{g,h,k})
= \check r_{g,h,k}^{ghk}(\tilde\g_{g,h,k})
\end{equation}
for some $\tilde \g_{g,h,k}\in i_{g,h,k}^{ghk}(A_{g,h,k})$.
So vice--versa to show associativity one needs to show that
\begin{equation}
\check r_{gh,k}^{ghk}
(r_{gh,k}^{gh}(\check r^{gh}_{g,h}(\tilde \g_{g,h}))\tilde \g_{gh,k})= 
  \check r_{g,h,k}^{ghk}(\tilde \g_{g,h,k})
\end{equation}
for some $\tilde \g_{g,h,k}$
which is a symmetric expression in the indices.

\subsection{Intersection $G$ Frobenius algebras}
Vice--versa in the given $G$--intersection Frobenius datum
 using the diagram (\ref{multdiagram}) 
as an Ansatz for a multiplication  we will arrive at a special
type of Frobenius algebra. The associativity of this Ansatz can then
be checked on the triple intersections.

\subsubsection{Definition}
An intersection $G$--Frobenius algebra is 
an intersection $G$--Frobenius datum of level $k\geq 3$
together with a $G$--Frobenius algebra structure on $A:= \bigoplus A_g$
whose multiplication is given by the diagram (\ref{multdiagram}) and
whose associativity is given by diagram (\ref{assdiagram})

\subsubsection{Remark}
\label{specialintersect} Reconstructing from special reconstruction data
one can define the algebras $A_{g_1,\dots g_n}$ via the 
following procedure. Set $I_{g_1,\dots, g_n}:= I_{g_1} + \dots +I_{g_n}$
and $A_{g_1,\dots g_n}:= A_e/I_{g_1,\dots, g_n}$. In order to
get $G$--intersection Frobenius data one has then only to show that the
 $A_{g_1,\dots g_n}$ are indeed Frobenius algebras and choose a metric
for them. If this is possible then Proposition \ref{intersect}
shows that any reconstructed special $G$ Frobenius algebra is 
an intersection $G$ Frobenius algebra.

\subsubsection{Examples} 
\begin{itemize}
\item[i)] We will show that the structures of Remark \ref{specialintersect}
are indeed present in the case of symmetric products.
\item[ii)] The $G$--Frobenius structures for the global orbifold cohomology
ring as presented in \cite{FG} are intersection $G$--Frobenius algebras.
\end{itemize}

\subsection{The Sign}
Given a preferred choice of character, it is possible to define a sign which
corresponds to a super--twist from a preferred choice of super--grading.

\subsubsection{Remark} 
Given a special $G$--Frobenius algebra $A$  
we denote the Eigenvalue of $\rho$ w.r.t.\ $\varphi_g$ by $\l_g$
and furthermore denote the Eigenvalue of $\varphi^h_g$ on $i_h(\rho_h)$ by
$\l_{g}^h$ i.e.\ $\varphi_g(\rho) = \l_g \rho$ and 
$\varphi^h_g(i_h(\rho_h))= \l_g^h i_h(\rho_h)$.
By the projective $G$--invariance of the metric 
\begin{equation}
\l_h = \chi_h^{-2} 
\end{equation}
and we can regard the ensembles $\l_g$ and $\l^g_h$ as characters.

\subsubsection{Definition}
We define  a sign $\sign$ to be an element of $\Hom(G,k^*)$.
Fixing an element $\sign  \in \Hom (G,\Z2)$ we
can define the associated character $\psi$ by
\begin{equation}
\psi(g):= (-1)^{\sign(g)} \chi_g 
\end{equation}

Vice--versa given a character $\psi \in \Hom(G,k^*)$ with the property that
$\psi^2=\chi^2$ we define the {\em sign given by $\psi$} to be 
\begin{equation}
(-1)^{\sign(g)}:= \chi_g \psi(g)^{-1}
\end{equation}

Finally, any choice of root of $\l$ defines a sign.

Given $\sign$ and $\sign^g$ for $A$ and $A^g$ for all 
$g,h \in G, [g,h]=e$ we set 
\begin{equation}
\nu(g,h) \equiv \sign(g) +sign^g(h) + \tilde h^g +\tilde g \; (2)
\end{equation}
$\sign$ and $\sign^g$
are said to be compatible if for all $h \in g$ 
\begin{equation}
\label{nucompat}
\nu(g,h)=\nu(gh,h)=\nu(h,g)=\nu(g^{-1},h)
\end{equation}

\subsection{Algebraic Discrete Torsion}
\label{dtors}
In certain situations it is also possible to
 distinguish one $G$--Frobenius algebra 
as initial under the action of discrete torsion. This is the case for instance
for Jacobian Frobenius algebras. 
In general, we can define a similar structure for intersection Frobenius
algebras, which then incorporates the trace axiom into the definition 
of discrete torsion. This shows that the compatibility with the
trace axiom in principle fixes the action up to a twist by discrete torsion.

Denote the centralizer of an element $g\in G$ by $Z(g)$ and fix
a sign of A. We will consider $G$--intersection Frobenius data
of level 2.

\subsubsection{The induced $Z(g)$--Frobenius algebra structure} If we are in an intersection Frobenius algebra
of level $k \geq 2$, 
given $A_g$ we can consider \\
{\sc The underlying additive structure.}\\ 
\begin{equation}
\hat A_g = \bigoplus_{h \in Z(g)} (A_g)_h = \bigoplus_{h \in Z(g)} A_{g,h}
\end{equation}
Notice that if $h\in Z(g)$, $\varphi_h: A_g \rightarrow A_g$ and 
$\varphi$ descends to a $Z(g)$ action on $A_g$. However,
we have that $\varphi_h (1_g) =  \varphi_{h,g}1_g$, but $1_g$ should be 
invariant under the $Z(g)$--action as the new identity.
Therefore we set

{\sc The $Z(g)$--action.}\\
\begin{equation}
\varphi_h^g:=\varphi_{h,g}^{-1}\varphi_h
\end{equation}
With this 
definition $\varphi^g_h (1_g) =\varphi_{h,g}^{-1}\varphi_{h,g}1_g =1_g$.

{\sc The character.}
Given a $G$--action on the level 2 $G$--intersection algebra,
we can augment the picture with a character $\chi_h^g$, which will be
determined by the trace axiom. 

{\sc Supergrading.}
We fix the super--degree of $A_{g,h}$ in  $\hat A_g$
and denote it by $\tilde h^g$.

\subsubsection{Definition} An intersection Frobenius algebra of level
$k\geq 2$ is said to
satisfy
{\em the discrete torsion condition}, if 
the above data satisfy the 
projective trace axiom and for all $g,h \in G$ there are isomorphisms between 
$A_{gh,h} \iso A_{g,h}$.

\subsubsection{Proposition} In an intersection Frobenius algebra $A$ of level
$k\geq 2$ that satisfies the discrete torsion condition, the following
equality holds for all $g,h \in G, [g,h]=e$:

\begin{equation}
\chi_g \STr(\varphi_g|_{A_h})= \varphi_{g,h}\chi_g (\chi_g^h)^{-1} (-1)^{\tilde g}
(-1)^{\tilde h^g}\dim(A_{g,h})
\end{equation}
or given roots $\psi,\psi^g$ of $\l,\l^g$
\begin{equation}
\label{dttrace}
\chi_g \STr(\varphi_g|_{A_h})=
\varphi_{g,h}\psi_g (\phi_g^h)^{-1}(-1)^{\sign(g)+\sign^h(g)} (-1)^{\tilde g}
(-1)^{\tilde h^g}\dim(A_{g,h})
\end{equation}

{\bf Proof.}
From the discrete torsion condition we obtain
$$
(-1)^{\tilde g^h} \dim(A_{g,h})= \chi_h^g \STr(\varphi^h_g|_{A_{h,e}})
$$
and furthermore
$$
\STr (\varphi_g|_{A_h}) = (-1)^{\tilde g} \varphi_{g,h}^{-1} 
\STr(\varphi^h_g|_{A_{h,e}})
$$

\subsubsection{Corollary} If $\psi$ and $\psi^g$ are compatible 
then
\begin{equation}
\chi_g \STr(\varphi_g|_{A_h})=
\varphi_{g,h}\psi_g (\phi_g^h)^{-1}(-1)^{\sign(g)+\sign(h)} 
(-1)^{\nu(g,h)}\dim(A_{g,h})
\end{equation}

\subsubsection{Definition} If $\sign$ and the $\sign^g$ are compatible, 
we set for $g,h\in G, [g,h]=e$
\begin{equation}
 T(h,g)
=(-1)^{\sign(g)\sign(h)} (-1)^{\sign(g)+\sign(h)}
(-1)^{|\nu_{g,h}|}\dim(A_{g,h})
\end{equation}
it satisfies for $g,h \in G, [g,h]=e$
\begin{equation}
T(g,h)=T(h,g)=T(gh,h)=T(g^{-1},h)
\end{equation}

\begin{equation}
\label{dt}
\eps(h,g)=\varphi_{g,h} (-1)^{\sign(g)\sign(h)}\psi_g (\psi_g^h)^{-1}
\end{equation}

Due to the projective trace axiom and by definition
$\eps$ viewed as a function from $G\times G \rightarrow k^*$
satisfies the conditions
of discrete torsion which are defined by:
\begin{equation}
\label{dtorsion}
\eps(g,h)=\eps(h^{-1},g)
 \quad \eps(g,g)=1 \quad \eps(g_1g_2,h)= \eps(g_1,h)\eps(g_2,h)
\end{equation}

\section{Jacobian Frobenius Algebras}

We first recall the main definitions and statements about Jacobian 
Frobenius algebras from [K2, K3].

\subsection{Reminder} A Frobenius algebra $A$ is called {\it Jacobian}
if it can be represented as the Milnor ring of a function $f$. I.e. if 
there is a function $f\in {\Cal O}_{{\bf A}^{n}_{k}}$ 
s.t. $A = {\Cal O}_{{\bf A}^{n}_{k}}/J_{f}$ 
 where $J_{f}$ is the Jacobian 
ideal of $f$. And the bilinear form is given by the residue pairing.
This is the form given by the Hessian of $\rho={\rm Hess}_{f}$.

If we write ${\Cal O}_{{\bf A}^{n}_{k}}= k[x_{1}\dots x_{n}]$, 
$J_{f}$  is  the ideal spanned by the $\frac{\del f}{\del x_{i}}$.

A {\it realization of a Jacobian Frobenius algebra} is a pair $(A,f)$
of a Jacobian Frobenius algebra and a function $f$  on some 
affine $k$ space ${\bf A}_{k}^{n}$, i.e. $f \in {\Cal O}_{{\bf 
A}_{k}^{n}}= k[x_{1}\dots x_{n}]$ s.t. $A= k[x_{1}\dots x_{n}]$ and
$\rho := {\rm det}(\frac{\del^{2} f}{\del x_{i}\del x_{j}})$.

\subsection{Definition} A {\it natural $G$ action on a realization of a
Jacobian Frobenius algebra $(A_e,f)$ } 
is a linear $G$ action on ${\bf A}_{k}^{n}$ which  leaves $f$ 
invariant. 
Given a natural $G$ action on a realization of a
Jacobian Frobenius algebra $(A,f)$ set for each $g\in G$,
${\Cal O}_{g}:= {\Cal O}_{ {\rm Fix}_{g}({\bf A}_{k}^{n})}$.

We also write $V(g):= {\rm Fix}_{g}({\bf A}_{k}^{n})$.

This is the ring of functions of the fixed point set
of $g$ for the $G$ action on ${\bf A}_{k}^{n}$. 
These are the functions fixed by $g$:
${\Cal O}_{g}= k[x_{1},\dots, x_{n}]^{g}$.

Denote by $J_{g}:= J_{f|_{{\rm Fix}_{g}({\bf A}_{k}^{n})}}$ the 
Jacobian ideal of $f$ restricted to the fixed point set of $g$.

Define
\begin{equation}
A_{g}:= {\Cal O}_{g}/J_{g}
\end{equation}
The $A_{g}$ will be called twisted sectors for $g \neq 1$.
Notice that each $A_{g}$ is a Jacobian Frobenius algebra with the 
natural realization given by $(A_{g}, f|_{\mathrm {Fix}_{g}})$.
In particular, it comes equipped with an invariant bilinear form 
$\tilde\eta_{g}$
defined by the element $\mathrm {Hess}(f|_{\mathrm {Fix}_{g}})$.

For $g = 1$ the definition of $A_e$ is just the
realization of the original Frobenius algebra,
which we also call the untwisted sector.

Notice there is a restriction morphism $r_{g}: A_e \rightarrow A_{g}$
given by $a \mapsto a \;{\rm mod} \;J_{g}$.

Denote $r_{g}(1)$ by $1_{g}$. This is a non--zero element of 
$A_{g}$ since the action was linear. 
Furthermore it generates $A_{g}$ as a cyclic $A_e$ module.

The set ${\rm Fix}_{g}{\bf A}_{k}^{n}$ is a linear subspace. 
Let $I_{g}$ be the vanishing ideal of this space.

We obtain a sequence
$$
0 \rightarrow I_{g} \rightarrow A_e\stackrel {r_{g}}{\rightarrow} 
A_{g}\rightarrow 0
$$

Let $i_{a}$ be any splitting of this sequence induced by the inclusion:
$\hat i_{g}: {\Cal O}_{g} \rightarrow {\Cal O}_{e}$ which descends due 
to the invariance of $f$.

In coordinates, we have the following description. 
Let ${\rm Fix}_{g}{\bf A}_{k}^{n}$ be
given by equations $x_{i}=0: i \in N_{g}$ for some index set $N_{g}$.

Choosing complementary generators $x_{j}: j \in T_{g}$, we
have ${\Cal O}_{g}= k[x_{j}:j \in T_{g}]$
and 
${\Cal O}_{e}= k[x_{j},x_{i}:j \in T_{g}, i\in N_{g}]$.
Then $I_{g}=(x_{i}: i \in N_{g})_{{\Cal O}_{e}}$ is 
the  ideal in ${\Cal O}_{e}$
generated by
the $x_{i}$ and
${\Cal O}_{e} = I_{g}\oplus  i_{g}(A_{g})$ using the splitting 
$i_{g}$ coming from the natural inclusion
$\hat i_{g}:k[x_{j}:j \in T_{g}] \rightarrow k[x_{j},x_{i}:j \in T_{g}, i\in N_{g}]$. 
We also define the projections 
$$
\pi_{g}: A_e \rightarrow A_{e}; \pi_{g} = i_{g} \circ r_{g}
$$ 
which in coordinates are given by $f \mapsto f|_{x_{j}=0: j \in 
N_{g}}$
Let 
$$
A := \bigoplus_{g \in G} A_{g}
$$

where the sum is a sum of $A_e$ modules.

Some of the conditions of the reconstruction program are automatic
for Jacobian Frobenius algebras. The conditions and freedoms of choice of 
compatible data to the above
special reconstruction data are given by the following:

\subsection{Theorem (Reconstruction for Jacobian algebras)} 
Given a natural $G$ action on a realization of a Jacobian Frobenius 
algebra $(A_e,f)$ with a quasi--homogeneous function $f$  
with $d_g= 0$ iff $g=e$
together with a natural 
choice of splittings $i_{g}$ the
possible structures of naturally graded 
special $G$ twisted Frobenius algebra on the 
$A_e$ module $A := \bigoplus_{g \in G} A_{g}$ are in 1--1
correspondence with the set of 
section independent $G$ graded cocycles $\g$ which are compatible with 
the metric together with a 
choice of sign $\sign\in Hom(G,\Z2)$ and
a compatible non--abelian two cocycle $\varphi$ with values in 
$k^{*}$, which satisfy the condition of discrete torsion

\begin{equation}
\label{dtcondition}
\forall g,h \text{ s.t.\ } [g,h]=e: 
\varphi_{g,h}\varphi_{h,g}\det (g|_{N_{h}})\det (h|_{N_{g}})=1
\end{equation}
\medskip
and the supergrading condition
\begin{equation}
|N_g|+|N_h|\equiv |N_{gh}| \; (2) \text{ or } \g_{g,h}=0
\end{equation}

This means in particular that the trace condition is replaced by 
(\ref{dtcondition}).
Also notice that if $\g_{g,h} \neq 0$ then the factor
$\varphi_{g,h}\varphi_{h,g}=1$ in (\ref{dtcondition})
by the compatibility equations so that (\ref{dtcondition}) reads
\begin{equation}
\det (g|_{N_{h}})\det (h|_{N_{g}})=1
\end{equation}

{\sc Notation.} If $[g,h]\neq 0$ then $\deg(g|_{N_{h}})$ is taken as
an abbreviation for\\
 $\deg(g) \det^{-1}(g|_{T_h})$.

\subsubsection{Character and Sign} The character and parity are
fixed by a choice of sign $\sign$ and are given by:

\begin{equation}
\chi_{g}= (-1)^{\tilde g} (-1)^{|N_{g}|} \mathrm{det}(g)
\end{equation}

The sign is defined by 

$$
\chi_{g}= (-1)^{\sign(g)} \det(g)
$$
i.e.\ we choose $\psi_g = \det(g)$

and satisfies
\begin{equation}
\label{jacsign}
\sign(g):= \tilde g +|N_g| \mod 2
\end{equation}

\subsubsection{Bilinear form on the twisted sectors}
If the character $\chi$ in non--trivial, we have to shift the
natural bilinear forms $\eta_g$ on $A_g$ by
\begin{equation}
((-1)^{\tilde g}\chi_g)^{1/2}\eta_g
\end{equation}
where we choose to cut the plane along the negative real axis.
For more comments on this procedure see [K3] and the following remarks.

\subsubsection{Remarks about the normalization}
\label{normremarks}
We would like to point out that the setup of reconstruction data
already includes the forms $\eta_g$. This is the reason for the above shift.
Indeed there is always a pencil of metrics for any given irreducible Frobenius 
algebra. The overall normalization is fixed by $\g_{g,g^{-1}}$. 
More precisely, we always have the equation:
\begin{equation}
\g_{g,g^{-1}}i_g(\rho_g)= \rho
\end{equation}
Notice that since $\g_{g,g^{-1}} I_g=0$ this equation determines
$\rho_g$ uniquely at least in the graded irreducible case since $\rho_g$ is of
necessarily of top degree in $A_g$.
So if, we were not to include the $\eta_g$ into the data, the only conditions on
 the $\g_{g,g^{-1}}$ would be that they do not vanish, live in the right degree and satisfy the compatibility but there would be  no need for rescaling.

Another way to avoid the shift is to include it in the restriction data
by setting

\begin{equation}
A_g:= \mathcal{O}_{f_g} \text{ with } f_g= ((-1)^{\tilde g}\chi_g)^{1/2}\eta_g f|_{\mathrm{Fix}(g)}
\end{equation}

\subsubsection{Natural discrete Torsion for Jacobian Frobenius algebras}

We can write 
$$
\chi_h {\rm STr} (\varphi_{h}|_{A_{g}})=\eps(h,g)  T(h,g)
$$
where
\begin{eqnarray}
 T(h,g)&=& (-1)^{\sign(g)\sign(h)} (-1)^{\sign(g)+\sign(h)} 
(-1)^{|T_{g}\cap T_{h}|  +N} \nn\\
&&\quad \dim(i_g(A_{g})\cap i_h(A_{h}))\nn\\
&=&(-1)^{\sign(g)\sign(h)} (-1)^{\sign(g)+\sign(h)}
(-1)^{|N_{g,h}|}\dim(A_{g,h})
\label{Tgh}
\end{eqnarray}
where we introduced the notation $|N_{g,h}|$ for 
$\dim(\mathrm{Fix}(g)\cap\mathrm{Fix}(h))$
and $A_{g,h}$ \\
for $\mathcal {O}_{f|_{\mathrm{Fix}(g)\cap\mathrm{Fix}(h)}}$
\begin{equation}
\label{jacdt}
\eps(h,g)=\varphi_{g,h} (-1)^{\sign(g)\sign(h)}\det (g|_{N_{h}})
\end{equation}
The projective trace axiom is satisfied in the graded case if $\eps$ satisfies 
the equations of discrete torsion

\begin{equation}
%\label{dtorsion}
\eps(g,h)=\eps(h^{-1},g)
 \quad \eps(g,g)=1 \quad \eps(g_1g_2,h)= \eps(g_1,h)\eps(g_2,h)
\end{equation}
which in terms of the $\varphi$ is equivalent to the condition 
(\ref{dtcondition}).

\subsubsection{Remark}
This definition of discrete torsion agrees with the more general one 
of \ref{dtors} if we set $\psi = \det(g)$ and $\psi^g(h)=\det(h)|_{T_g}$.
Indeed we find $sign^g(h) \equiv  \tilde h^g +|N^g_{g,h}|$ 
with  $|N^g_{g,h}|= \codim_{\Fix_g}(\Fix_g \cap \Fix_h)$
and thus
\begin{eqnarray}\nu(g,h)&\equiv& 
\sign(g) +sign^g(h) + \tilde h^g +\tilde g \; (2) \nn\\
&\equiv& \sign(g) +\codim(\Fix_g \cap \Fix_h)+ |N_g|+\tilde g
\equiv |N_{g,h}|\; (2)
\end{eqnarray}
 
\subsubsection{Examples}
\begin{itemize} 
\item[1)]\label{point}  (pt/G). Recall (cf.\ [K3]) 
that given a linear representation
$\rho: G \rightarrow O(n,k)$, we obtain the G--twisted Frobenius algebra
$pt/G$ from the Morse function $f=z_{n}^{1}+ \ldots + z^{2}_{n}$.

All sectors are isomorphic to $k$:

$$A= \bigoplus_{g\in G} k$$ 
all the $d_{g}=0$
and all the $r_{g}=id$. In particular, we have that $\g_{g,g^{-1}}= 
\check r_{g}(1) =1$ and $\pi_{g}(\g_{h,h^{-1}})=1 \neq 0$, so we see 
that the $\g_{g,h}\in k^{*}$ and are given (up to rescaling) by group 
cocycles $\g \in H^{2}(G,k^{*})$ and since the $g_{g,h}\neq 0$ the
$\varphi$ and hence the discrete torsion are fixed by the compatibility
$\g_{g,h}= \varphi_{g,h}\g_{ghg^{-1},g}$.

Explicitly: Fix a parity $\; \tilde {} \in \Hom(G,\Z2)$.

The sign and character are given by 
\begin{equation}
\sign(g) \equiv \tilde g  \quad \chi_g = (-1)^{\sign(g)}=  (-1)^{\tilde g}
\end{equation}

\item[2)] Another example to keep in mind is $A_{n}$ which is the
Frobenius algebra associated to $z^{n+1}$ together with the $\Znn$ action
$z \mapsto \zeta_n z$ where $\zeta_n^n=1$ [cf.\ K3].

\item[3)] $A^{\otimes n}$ together with the permutation action. We will consider this example in depth in \S \ref{snjac} and \S \ref {symp}.
This example has appeared many times in different guises in 
[DHVV,D1,D2,LS,U,WZ]. 
Our treatment is the completely general and subsumes all these cases. Also,
there is an ambiguity of
signs which is explained by our treatment.
\end{itemize}

\subsection{Theorem} Jacobian algebras naturally give intersection algebras.

{\bf Proof.} This is straight--forward. We set 
\begin{equation}
A_{g_1, \dots, g_k} := \mathcal{O}_{f_{g_1, \dots, g_k}} \text { with }
 f_{g_1, \dots, g_k} := f|_{\bigcap_{i=1}^k\mathrm{Fix(g_i)}}
\end{equation}
and use the obvious restriction maps. Here again the remarks of 
\ref{normremarks} apply.

\section{Special $\Sn$--twisted Frobenius algebras}
\label{sn}

\subsection{Notation}
\label{snnotation}
Given a permutation $\s \in \Sn$, we associate to it its cycle decomposition
$c(\s)$ and its index type $I(\s) := \{I_1, \dots I_k\}$ where the $I_j$
 are the independent 
sets in the cycle decomposition of $\s$.
Notice that the $I(\s)$ can also be written as $\la \s \ra \backslash \bar n$
where this is the quotient set of $\bar n$ w.r.t.\ group action of the 
group generated by $\s$.

The length of a cycle decomposition 
$|c(\s)|$ is defined to be the number of independent cycles in the 
decomposition. The partition gives rise to its norm $(n_1, \dots, n_k)$ of $n$
where $n_i :=|I_i|$. And the type of a cycle is defined to be
$(N_1(\s),N_2(\s), \dots )$ 
where $N_i= \#\text{ of } n_j =i$ in $(n_1, \dots, n_k)$,
i.e.\ $N_i$ the number of cycles of length $i$ in the cocycle decomposition of
$\s$.

 We define the degree of $\s \in \Sn$ to be
$|\s| :=$  the minimal length of $\s$ as a word in transpositions 
$= n-|c(\s)|$. 

Recall the relations in $\Sn$ are

\begin{eqnarray}
\tau^{2}&=&1\\
\t\t'&=&\t'\t'' \mbox{ where } \t = (ij), \t'=(jk),\t''=(kl).
\label{relation}
\end{eqnarray}

\subsection{Definition} We call two elements 
$\sigma,\sigma' \in \Sn$ {\em transversal}, if 
$|\sigma\sigma'|=|\sigma|+|\sigma'|$.

\subsection{The linear subspace arrangement}
\label{vsigmas}
A good deal of the theory of $\Sn$ Frobenius algebras is governed
by the canonical permutation representation of $\Sn$ on $k^n$ given
by $\rho(\s)(e_i)= e_{\s(i)}$ for the canonical basis $(e_i)$ of $k^n$.

We set $V_{\s}:=\mathrm{Fix}(\s)$

and $V_{\s_1, \dots, \s_n} := \bigcap_{i=1}^n V_{\s_i}$.
Notice that
\begin{equation}
l(\s) =  \dim( V_{\s})= |\la \s \ra\backslash \bar n| 
\end{equation}
and 
\begin{equation}
|\s| = \codim(V_{\s})
\end{equation}
In the same spirit, we define

\begin{eqnarray}
l(\s_1, \dots \s_n)&:=& \dim(V_{\s_1, \dots, \s_n})\nn\\
|\s_1, \dots, \s_n|&:=&\codim(V_{\s_1, \dots, \s_n})
\end{eqnarray}

This explains the name  transversal. 
Since if $\s$ and $\s'$ are transversal then

$$V_{\s,\s'}=V_{\s}\cap V_{\s'}=V_{\s\s'}$$ 
and the intersection is transversal.

Furthermore notice that 
\begin{equation}
l(\s_1, \dots, \s_n) = |\la \s_1, \dots, \s_n \ra \backslash \bar n|.
\end{equation}
where again the last set is the quotient set of $\bar n$ 
by the action under the 
group generated by $\s_1, \dots,\s_n$.

\subsection{Definition}
We call a cocycle  $\g:\Sn \times \Sn \rightarrow A$ 
{\em normalizable} 
if for all {\em transversal} pairs  $\t, \s \in \Sn , |\t|=1:
\g_{\s,\t}\in A_e^*$, i.e.\ is $\g_{\s,\t}$ is invertible, 
and {\em normalized} if it is normalizable and
 for all {\em transversal } $\tau, \sigma \in \Sn , |\t|=1:
\g_{\s,\t}=1$.

In the example of symmetric products of an irreducible Frobenius
algebra or in general $A_e$ irreducible the invertibles are of 
degree 0 and are given  precisely by $k^*$.

\subsubsection{Lemma} 
\label{transnorm}
If a cocycle is normalized then for any {\it 
transversal} 
$\sigma,\sigma' \in \Sn: \g_{\sigma, \sigma'}=1$.

{\bf Proof.}
We write $\sigma'= \t'_{1}\cdots \t'_{k}$ 
with $k=|\sigma'|$
where all $\tau_{i}$
are  transpositions.

Thus by associativity: 
$$\s\s'= (((\dots(\s\tau'_{1})\tau'_{2}) \cdots )\tau'_{k}),$$ 
so
$$\g_{\s\s'}=\pi_{\s\s'}(\g_{\s,\s'}) = 
\pi_{\s\s'}(\prod_{i=1}^{k}
\g_{\s\prod_{j=1}^{i-1}(\tau_{j}),\t_i}) = 
\pi_{\s\s'}(\prod_{i=1}^{k} 1) =1$$

\subsubsection{Remark} 
Recall that  $\g_{\t,\t}= \check{r}_{g}(1_{\t})$ for a
transposition
$\t$.

\subsubsection{Lemma} Let $\s \in \Sn$. If $\g$ is a normalized 
cocycle, then
for any decomposition into transpositions $\s=\tau_{1}\cdots 
\tau_{|\s|}: \g_{\s,\s^{-1}}\prod_{i=1}^{|\s|}\g_{\t_{i},\t_{i}}$

{\bf Proof.} Let $k=|\s|$.
Thus by associativity: 
$$\s\s^{-1}= (\tau_{1}(\tau_{2}( \cdots 
(\tau_{k}\tau_{k}\cdots\tau_{2}\tau_{1})\cdots ))),$$
and if $\t$ and $\s'$ are transversal
$$\pi_{\s'}(\g_{\t,\t\sigma'})=
\pi_{\s'}(\g_{\t,\t\s'}\g_{\t,\s'})= \pi_{\s'}(\g_{\t,\t}\g_{e,\s'})
=\pi_{\s'}(\g_{\t,\t}).$$
So $\g_{\s,\s^{-1}}= \prod_{i=1}^{k} \g_{\tau_{i},\tau_{i}}$.

\subsection{Theorem}
\label{unique} Given special $S_{n}$ reconstruction data,
a choice of normalized cocycle $\g:\Sn \times \Sn \rightarrow A$ is
unique. Furthermore a choice of normalizable cocycle is fixed by a
choice of the $\g_{\t,\s}$ with $\t$ and $\s$ transversal.

{\bf Proof.}
We have that the $\g_{\s,\s^{-1}}$ are given by 
$\g_{\s,\s^{-1}}= \check{r}_{\s}(1_{\s})$ and thus fixed after the 
normalization which fixes the $r_{\s}$.
Again choosing any minimal decomposition
$\s'=\tau'_{1}\cdots 
\tau'_{|\s|}$
and by using the normalization and associativity repeatedly, we obtain that
\begin{eqnarray*}
\g_{\s,\s'}=\pi_{{\s\s'}}(\g_{\s,\s'}\prod_{i=1}^{|\s'|}
\g_{\tau'_{i+1},\prod_{j=1}^{i}\tau'_{j}})
=\pi_{\s,\s'}(\prod _{i=1}^{|\s'|}
\g_{\s\prod_{j=1}^{i}\tau'_{i-1},\tau'_{i}})
\\
=\pi_{\s\s'}(\prod_{i \in I}\g_{\tau'_{i},\tau'_{i}}) 
\end{eqnarray*}
where  
$I:=\{i: |\s(\prod_{j=1}^{i-1}\tau'_{j})\tau'_{i}|
=|\s\prod_{j=1}^{i-1}\tau'_{j}|-2\}$.

Thereby the $\g_{\s,\s'}$ are already determined by the $\g_{\t,\t}$
which are in turn given by $\check{r}_{\tau}(1_{\tau})$.

If the cocycles are only normalizable, we obtain the result in a similar
fashion.

\subsection{Discrete torsion for $\Sn$}

It is well known (see e.g.\ \cite{Ka}) that $H^2(\Sn,k^*)=\Z2$.

\subsection{Lemma}
\label{sntorsion}
 Let $\Phi$ be a cocycle corresponding
to the non--trivial central extension  of $\Sn$ defined as
the group generated by $\hat \t_i:i=1,\dots n$
$$
\hat\t_i\hat\t_i=z, \quad zz=e,\quad  
\hat\t_i\hat\t_{i+1}\hat\t_i =  \hat\t_{i+1} \hat\t_i\hat\t_{i+1},\quad 
\hat\t_i\hat\t_j = z \hat\t_j\hat\t_i: |i-j|\geq2
$$
and let $k^{\Phi}[\Sn]$ be the corresponding twisted group ring 
(here $z\mapsto -1$) 
then
$$
\eps_{\Phi}(\t_i,\t_i)=1 \quad \eps(\t_i,\t_j)= -1 : i\neq j
$$

{\bf Proof} 
Since $\hat\t_i^2=-1$, $\hat\t_i\hat \t_i \hat\t_i^{-1}=
-(-\hat \tau_i)=\hat \tau_i$. If $|i-j|\geq 2$
$\hat \tau_i \hat\tau_j\hat\tau_i^{-1}=- \hat\tau_j\hat \tau_i(-\hat\tau_i)
= - \hat\tau_j$.

\subsection{Supergrading and Parity $p$}
\label{snsuper}
Since $\Sn$ is generated by transpositions which all lie in the 
same conjugacy class, we see that the choices of $\Z2$--grading 
$\; \tilde{ } \in \Hom(\Sn, \Z2)$ are 
given by 
\begin{itemize}
\item[i)] pure even $\forall \s: \tilde \s=1$. We call this the even case
and set the parity $p=0$.
\item[ii)] The sign representation $\tilde \s \equiv |\s|\; (2)$.
We call this the odd case
and set the parity $p=1$.
\end{itemize}

\subsection{Lemma} For the (super) twisted group ring,
the following equations hold:

$$
\eps(\s,\s')=(-1)^{p|\s||\s'|}
$$
in particular $\forall \t,\t'\in \Sn, |\t|=|\t'|=1, [\t,\t']=e$
$$
\eps(\t,\t)=(-1)^p\quad \eps(\t,\t')=(-1)^p
$$

This follows from the general result \ref{supergroupring}

\subsection{The non--abelian cocycles $\varphi$}
\subsubsection{Remark}
Due to the relation (\ref{nonabcocycle}),
we see that $\varphi$ is determined by the $\varphi_{\t,\s}$ 
with $|\t|=1$.

\subsubsection{Lemma} 
\label{pdef} For any non--abelian $\Sn$ cocycle 
$\varphi$ there is a fixed $p\in \{-1,+1\}$ s.t.\ for  all $\t\in \Sn,|\t|=1$
$\varphi_{\t,\t}=(-1)^{\tilde \t \tilde \t}=(-1)^p$. Furthermore if $\varphi$
is compatible with a section independent cocycle compatible with the metric,
then $p$ is the supergrading as an element in $\Zz$ (see \ref{snsuper}).

{\bf Proof.} By the definition of a non--abelian cocycle, we see that
$\forall \t: \varphi_{\t,\t}\in \{-1,1\}$. Furthermore
all transpositions are conjugate so that by
\ref{doubleconj} $\varphi_{\t,\t}=\varphi_{\t',\t'}$
for $\t,\t'\in \Sn:|\t|=|\t'|=1$ which shows the claim. In the case
of a compatible pair furthermore:
$\g_{\t,\t}=\varphi_{\t,\t}(-1)^{\tilde \t \tilde 
\t}\g_{\t,\t}$ and $\g_{\t,\t}\neq 0$, so that 
$\varphi_{\t,\t}=(-1)^p=(-1)^{\tilde \t}$.

\subsubsection{Lemma} 
\label{qdef}
For $\t,\t'\in \Sn, \t\neq\t', |\t|=|\t'|=1 , [\t,\t']=e$
$\varphi_{\t,\t'}= (-1)^q$
for a fixed $q\in \{-1,1\}$.

{\bf Proof.} Since $\t\t=[\t,\t']=e$, by  (\ref{nonabcocycle}) 
$\varphi_{\t,\t'}=\pm 1$ and by (\ref{doubleconj}), the
value is indeed fixed simultaneously for all commuting transpositions,
since all pairs of commuting transpositions  are conjugate to each other. 

\subsection{Definition} 
We call a non--abelian cocycle $\varphi$ 
normalizable if 
for all $\t,\t'\in \Sn, \t\neq\t', |\t|=|\t'|=1,[\t,\t']=e, \varphi_{\t,\t'}=(-1)^p$
for some fixed $p \in \{-1,1\}$.

We call a non--abelian cocycle $\varphi$ 
normalized if $\forall \s,\t\in \Sn, |\t|=1$

$\varphi_{\s,\t}= (-1)^{\tilde \s\tilde \t}= (-1)^{p|\s|}$.

\subsection{Lemma}
After a possible twist by any  discrete torsion $\a$ with $[\a]\neq 0$   all 
non--abelian cocycle $\varphi$ 
normalizable. 

{\bf Proof.} By  Lemmas \ref{pdef} and \ref{qdef}, we have that indeed
for $\t,\t'\in \Sn, |\t|=|\t'|=1,[\t,\t']=e$
$\varphi(\t,\t) =(-1)^p$ and $\varphi(\t,\t')=(-1)^q$ with
$p,q\in \{-1,1\}$. If $p=q$ then the cocycle $\varphi$ 
is already normalizable.
If $p\neq q$, let $\Phi\in Z^2(\Sn,k^*)$ be  
the class given in  Lemma 
\ref{sntorsion} then $\varphi^{\Phi}(\t,\t) =(-1)^p$ and
$\varphi^{\Phi}=(-1)^p$ since if $p\neq q$ then $p=q+1$.
But on commuting elements $\eps_{\a}$ only depends on the
cohomology class of $\a$ and thus we could use a twist by
$\a$ for 
any class with
$[\a]\neq 0 \in H^2(\Sn,k^*)$ instead of $\Phi$.

If $\varphi$ is the non--abelian cocycle
of a special $\Sn$ Frobenius algebra $A$ then the non--abelian
cocycle $\varphi^{\Phi}$ can be obtained via tensoring 
with $k^{\Phi}[\Sn]$ as the non--abelian cocycle of
 $A^{\Phi}$.

The Theorem \ref{normal} 
contained in the Appendix A implies that all normalizable 
 non--abelian cocycles $\varphi$
can be rescaled to a normalized cocycle.

\subsection{Theorem}
\label{normalize} 
Any normalizable graded $\Sn$ cocycle $\g$ with
normalized
$\varphi$
can be normalized by a 
rescaling $1_{\sigma} \mapsto \lambda_{\sigma}1_{\sigma}$. 

And vice--versa given 
any normalized $\Sn$ cocycle and a choice of parity $p\in \{0,1\}$ there is only one compatible 
non--abelian cocycle $\varphi$ given by
\begin{equation}
\varphi_{\s,\s'}= (-1)^{p|\s||\s'|}
\end{equation}

{\bf Proof.}
First notice that by assumption of normalizability the 
$\g_{\s,\t}\in k^*$ for transversal $\t,\s$ 
we define the rescaling inductively on $|\s|$
by $\l_{\t}:=1$ and $\l_{\s} :=\l_{\s'}\g_{\s',\t'}$ where
$\s=\s'\t'$ and
$\t$ and $\s$ are transversal.

More precisely: 
let $\s= \s'\t'$ with $|\t|=1,|\s'|=|\s|-1$. With induction on $|\s|$
we define 
\begin{equation}
\label{scaling}
\lambda_{\s}:=\l_{\s'}\g_{\s',\t'}
\end{equation}

Then after scaling we obtain: 
$$\tilde \g_{\s',\t'} 
= \frac{\l_{\t'}\l_{\s'}}{\l_{\s'\t'}}\g_{\s',\t'}=\l_{\t'}=1$$

We have to show that (\ref{scaling}) is well defined i.e.\ is 
independent of the decomposition.
This  can again be seen by induction. 

First notice that if $|\s|= 1$, $\l_{\s}=1$ poses no problems.
If $|\s|=2$ either there is a unique decomposition into two disjoint
transpositions
or 
\begin{equation}
\label{rel}
\s= \t\t'=\t'\t''
\end{equation}
where $\t = (ij), \t'=(jk),\t''=(kl)$.
The first case again poses no problem. For the second one notice that
$\l_{\t}= \l_{\t'}=1$ and $\t'\t''\t'=\t$ thus

\begin{equation}
\label{transpose}
\g_{\t',\t''} = \varphi_{\t',\t''}(-1)^{\tilde \t' \tilde\t''}\g_{\t'\t''\t',\t'}= 
 \varphi_{\t',\t''}(-1)^{\tilde \t' \tilde\t''}\g_{\t,\t'}=\g_{\t,\t'}
\end{equation}

Assume the $\l_{\s}$ are well
defined for $|\s| < k$. Fix $\s$ with $|\s|=k$ and decompose 
$\s= \s'\t' =\s''\t''$ in two different ways.
Then we have to show that 
$$\l_{\s'}\g_{{\s'},\t'} = 
\l_{\s'}\g_{{\s''},\t''}$$
where by induction  $\l_{\s'}= 
\prod_{i=1}^{|\s'|}\g_{\prod_{j=1}^{i-1}\t'_{j},\t'_{i}}$ and 
$\s'=\prod_{i=1}^{|\s'|}\t'_{i}$ is any minimal representation.
We observe that in $\Sn$ we can obtain $\s'\t'$ from $\s''\t''$ by using the 
relation (\ref{rel}) repeatedly.
Thus by using associativity and (\ref{transpose})  we obtain:
$$\l_{\s''}\g_{{\s''},\t''}= 
(\prod_{i=1}^{|\s''|}\g_{\prod_{j=1}^{i-1}\t''_{j},\t''_{i}}) \g_{{\s''},\t''}
=(\prod_{i=1}^{|\s'|}\g_{\prod_{j=1}^{i-1}\t'_{j},\t'_{i}}) \g_{{\s'},\t'}=
\l_{\s'}\g_{{\s'},\t'}$$

The fastidious 
reader can find the explicit case study in Appendix B.

For the second statement notice that by Lemma \ref{transnorm} given a normalized $\g$ we 
have for {\em all} transversal $\s,\s': \g_{\s,\s'}=1$.

Thus for transversal $\t,\s$
$$
1=\g_{\t,\s}=\varphi_{\t,\s}(-1)^{\tilde \t \tilde \s}\g_{\t\s\t,\t}
= \varphi_{\t,\s}(-1)^{\tilde \t \tilde \s}$$
since $\t\s\t$ and $\t$ are transversal $|\t\s\t|=|\s|,|\t\s\t\t| = 
|\t\s|=|\t|+|\s|$.

And if $\s,\t$ are not transversal, then $\s= \t\s'$ with 
$|\s'|=|\s-1|$ and $\s'$ and $\t$ transversal.
$$\g_{\t,\t}=\g_{\t,\s}=\varphi_{\t,\s} (-1)^{|\s|} \g_{\s'\t,\t}= 
(-1)^{|\s|}\g_{\t,\t}$$
and since $\g_{\t,\t}\neq 0$, 
we find
\begin{equation}
\label{taus}
\varphi_{\t,\s}=\varphi_{\s,\t}=(-1)^p.
\end{equation}

And finally if $\s=\prod_{i=1}^{|\s|}\t_i$
$$
\varphi_{\s,\s'}=\prod_{i=1}^{|\s|}\varphi_{\t_i,\tilde \s'_i}
= (-1)^{p|\s||\s'|}
$$
by using (\ref{taus}) with
$\tilde \s_i= (\prod_{j=i+1}^{|\s|}\t_i) \s'(\prod_{j=i+1}^{|\s|}\t_i)^{-1}$,
$|\tilde \s_i|=|\s'|$.

\section{Symmetric powers of Jacobian Frobenius algebras}
\label{snjac}

In this paragraph, we study $\Sn$ orbifolds of $A^{\otimes n}$
where $A$ is a Jacobian Frobenius algebra.
We also fix the degree $d$ of $A$ to be the degree of $\rho$ ---
the element defining $\eta$.

The most important result for Jacobian Frobenius algebras (or manifolds)
is that $A_f \otimes A_g = A_{f+g}$ [K1]. Therefore 
$$
A_{f(z)}^{\otimes n} = A_{f(z_1) + \dots + f(z_n)}
$$
where $z$ is actually a multi-variable $z= (z^1,\dots ,z^m)$.   

\subsection {Remark}
\label{tensnot}
In the above notation, we should keep it mind that for 
functions $g_1,\dots g_n$, we have that
$$
g_1 \otimes \dots \otimes g_n  =g_1(z_1) \cdots g_n(z_n)
$$

\subsection{$\Sn$--action}
In this situation there is a natural action $\rho$ of $\Sn$ by permuting
the $z_i$ i.e.\ for $\s \in \Sn$
$$
\rho(\s)(z_i^k) = z^k_{\s(i)}
$$
It is clear that the function $f_n := f(z_1) + \dots + f(z_n)$
is invariant under this action, so that we can apply the theory of [K2, K3].
We see that the representation $\rho$ is just 
the $\dim A$--fold sum
of the standard representation of $\Sn$ on $k^n$.

\subsection{The twisted sectors}
To analyze the twisted sectors, we have to diagonalize the 
given representation. To this end, we regard the cycle decomposition
and realize that for each cycle with index set $I_l$ 
there is a $m$--dimensional Eigenspace
generated by
 $$\frac{1}{n_i}\sum_{i \in I_l}  z_i^l \text { for }  l=1,\dots, m$$
The other Eigenvectors being given by 
$$\frac{1}{n_l}\sum_{i \in I_l} \zeta_{n_l}^jf(i)  z_i^l$$
with Eigenvalue $\zeta_{n_l}^j$ where $f:I_l \rightarrow \{1,\dots,n_l\}$
is a bijective map respecting the cycle order.

Restricting $f_n$ to the space where all the variables with Eigenvalue 
different from one vanish see that 
$$
f_{\sigma}= f(z_i =z_j= u _k) \text{ if }  i,j \in I_k
$$

Using the variables $u_k$ it is obvious that

$$ 
A_{\s} = A_{f_{\s}} \iso A^{\otimes |\s|}
$$

\subsection{Restriction maps}
With the above choice of $u_k$ as variables and using Remark \ref{tensnot},
we find that the
restriction maps are given as follows:

$$
r_{\s}(g_1 \otimes \dots \otimes g_n) := \bigotimes_{i=1}^k( \prod_{j\in I_i} g_i)
\in A^{\otimes |\s|}
$$ 

Thus these maps are just contractions by multiplication.

\subsection{Fixed point sets}
By the above, we see that 
\begin{equation}
\mathrm{Fix(\s)}= \bigoplus_{i=1}^m V_{\s} \subset (k^n)^m
\end{equation}
where we used the notation of \ref{vsigmas}.  
Notice that 

\begin{equation}
\dim(V_{\s}) = m l(\s) \quad \codim(V_{\s})=|N_{\s}|=m|\s|
\end{equation}

\subsection{Bilinear form on $A^{\otimes n}$}
We notice that if the bilinear form on $A$ is given by the element
$\rho= \mathrm {Hess}(f)$ then 
 the bilinear form on $A^{\otimes n}$
is given by $\rho^{\otimes n}= \mathrm {Hess}(f_n)$ and it is 
invariant under the $\Sn$ action.
Indeed $\det^2(\rho(\s))=1$.
To be more precise, we have that 

$$
\det(\rho(\s))=(-1)^{m|\s|}
$$

(Here $\rho$ is of course the representation, not the element defining
the bilinear form.)

\subsection{The Character and Sign}
Notice that the character is either the alternating or the trivial one
depending on the choice of the sign, which is determined
by the choice of parity $p$ and on the choice 
of the number of variables $m$. (We have to keep in mind that we can always 
stabilize the function $f$ by adding squares of new variables).

Using the equation (\ref{jacsign}), we find however:
\begin{equation}
\chi_{\s}=(-1)^{\tilde \s}(-1)^{m|\s|}\det(\s)=(-1)^{\tilde \s}
\end{equation}
and find the sign of $\s$
to be 
\begin{equation}
\sign(\s)\equiv \tilde{\s} +m|\s|=(m+p)|\s| 
\end{equation}

Thus only the sign, but not the character depends on the number of variables!

\subsection{Bilinear form on the twisted sectors}
Since it is always the case that $(-1)^{\tilde \s}\chi_{\s}=1$,
we do not have to shift the natural bilinear forms on the twisted
sectors. They are given by $\eta^{\otimes l(\s)}$ or equivalently by
$\rho_{\s}= \rho^{\otimes l(\s)}$.

\subsection{Remark}
Notice also that since $\det(\rho(\s))= \pm 1$ 
(i.e.\ the Schur--Frobenius indicator is 1)
the form $\eta$ will
descend to the $\Sn$ invariants (see  e.g.\ [K3]).

\subsection{Proposition}
\label{normalizable}
 After a possible twist by discrete torsion
any compatible cocycle $\g$ is normalizable.

{\bf Proof.} We  check that 
$\pi_{\sigma}(\g_{\t,\t}) \neq 0$ for $\tau$ and $\sigma$ transversal.
Then the claim follows from Proposition \ref{zerocheck}.

Suppose $\tau$ and $\sigma$ are transversal and say $\tau = (ij)$, then
$i$ and $j$ belong to different subsets 
of the partition $I(\s)$  (say $I(\s)_i$ and $I(\s)_j$).
So since $\g_{\t,\t}\neq 0$ neither is $\pi_{\s}(\g_{\t,\t})$.

More explicitly: 
\begin{equation}
\g_{\tau,\tau}=\check r_{\tau}(1_\tau) = \sum_k 1 \otimes \dots \otimes 1
\otimes \stackrel {\stackrel {i}{\downarrow}}{a_k} \otimes 1 \otimes \dots
\otimes 1 
\otimes \stackrel{\stackrel {j}{\downarrow}}{b_k} \otimes 1 \otimes \dots
\otimes 1 
\end{equation}
where $\sum_k a_k \otimes b_k = \Delta(1) \neq 0 \in A\otimes A$ and 
$\Delta := \check\mu: A\rightarrow A\otimes A$ is the natural
co--multiplication on $A$.
And  
\begin{equation}
r_{\s}(\g_{\t,\t})=
 \sum_k 1 \otimes \dots \otimes 1
\otimes \stackrel {\stackrel {I(\s)_i}{\downarrow}}{a_k} \otimes 1 \otimes \dots
\otimes 1 
\otimes \stackrel{\stackrel {I(\s)_j}{\downarrow}}{b_k} \otimes 1 \otimes \dots
\otimes 1 
\end{equation}

Thus $\g_{\t,\t}$ is not in the kernel of the contraction $r_{\s}$
and thus not in the kernel of $\pi_{\s}$.

\subsection{Algebraic discrete Torsion} 
The choices of algebraic discrete torsion are given by the choices of 
cocycles $\varphi$ and the sign.
Since there is only one $\varphi$ for a given choice of parity and fixing
the parity the sign is determined by the number of variables $m$.

Recall (\ref{jacdt})
$$\eps(\s,\s')=\varphi_{\s,\s'}(-1)^{\sign(\s)\sign(\s')}\det(\s|_{N_{\s'}})
=(-1)^{m|\s||\s'|}\det(\s|_{N_{\s'}})
$$
and
\begin{eqnarray*}
T(\s,\s') &=& (-1)^{\sign(\s)\sign(\s')}(-1)^{\sign(\s)+\sign(\s')}
(-1)^{m|\s,\s'|}
\dim(A_{\s,\s'})\\
&=& (-1)^{p(|\s|+|\s'|+|\s||\s'|)}(-1)^{m(|\s|+|\s'|+|\s||\s'|+|\s,\s'|)}
\dim(A_{\s,\s'})
\end{eqnarray*}

\subsection{Reminder}
\label{cent}

Recall that the centralizer of an element $\s \in \Sn$ is
given by 
$$Z(\s) \cong \prod_k {\bf S}_{N_k} \ltimes {\bf Z}/k{\bf Z}^{N_k}$$ 
where $N_i$ the number of cycles of length $i$ in the cycle decomposition 
of $\s$ (cf.\ \ref{snnotation})
This result can also be restated as: ``discrete torsion can be undone
by a choice of sign''.

We note that $Z(\s)$ is generated by elements of the type $\t_k$ and $c_k$
where  $\t_k$ permutes two cycles of length $k$ of $\s$ 
and  $c_k$ is a cycle of length $k$ of $\s$.

Also $\eps$ is a group homomorphism in both variables, so
that by \ref{dtorsion} $\eps$ is 
fixed by its value on elements of the above type.

\subsubsection{Proposition} The discrete torsion is given by
$$
\eps(\s',\s)=\begin{cases} (-1)^{mk|\s|}(-1)^{m(k-1)}&
\text{if $\s'=\t_k$} \\
(-1)^{m(k-1)(|\s|-1)}& 
\text{if $\s'=c_k$}
\end{cases}
$$
where $\t_k$ and $c_k$ are the generators of $Z(\s)$ described above.

{\bf Proof.} 
$$
\det(\t_k)|_{N_{\s}} = \det(\t_k) \mathrm{det}^{-1}(\t_k|_{T_{\s}})= (-1)^{mk}(-1)^m
$$
and
$$
\det(c_k)|_{N_{\s}} = \det(c_k) \mathrm{det}^{-1}(c_k|_{T_{\s}})= (-1)^{m(k-1)}
$$

\subsubsection{Remark}
What this calculation shows is that we are dealing with the $m$--th
power of the non--trivial cocycle which in
the case $m=1$ has been calculated in [D2].
We again  see  the phenomenon that
the  addition of variables (stabilization) changes the sign and hence
the discrete torsion --- as is well known in singularity theory.
Actually the whole trace i.e.\ the product of $\eps$ and $T$ is constantly
equal to $(-1)^{p(|\s||\s'|+|\s|+|\s'|)}\dim(A_{\s,\s'})$
which coincides with the general statement c.f.\ (\ref{tracevalue}).

\subsubsection{Corollary}
The discrete torsion condition holds.

\subsection{Grading and shifts}
\subsubsection{Proposition}
\begin{eqnarray}
s^+_{\s}&=& d|\s| , \quad s^-_{\s}= 0 \\
s_{\s} &=&\frac {1}{2}(s_{\s}^+ + s_{\s}^-)= \frac{d}{2}|\s|
\end{eqnarray}
where $s^+$ and $s^-$ are the standard shifts for Jacobian Frobenius
algebras as defined in [K2,K3].

For the calculation of $s^+$, 
we fix some $\s \in \Sn$. Let $c(\s)$ be its cycle decomposition
and $I(\s):= \{I_1, \dots I_k\}$ be its index decomposition.
Then the shift $s^+_{\s}$ can be read off from the definition and the
identification
$$A_{\s} \iso \bigotimes_{i = 1}^{|c(\s)|} A_{I_i} \iso A^{n-|\s|}$$
with the degree of $A^{\otimes l}$ being $dl$, we obtain
$$
s^+_{\s} =nd - (n-|\s|)d = d|\s|
$$

The shift $s^-_{\s}$ is again calculated via the natural representation
$\rho: \Sn \rightarrow GL(n,k)$.

Recall (cf.\ [K3])

\begin{multline*}
s_{g}^- := \frac{1}{2\pi i}\mathrm{Tr} (\log(g))-\mathrm{Tr}(\log(g^{-1})):= 
\frac{1}{2\pi i}(\sum_i \l_i(g)-\sum_i \l_i(g^{-1}))\\
=\sum_{i: \l_i \neq 0} 2 (\frac{1}{2\pi i}\l_i(g)-1)
\end{multline*}

For a cycle $c$ of length $k$, we have the eigenvalues 
$\zeta_k^i,i=0, \dots k-1$ where $\zeta_k$ is the k-th root of unity
$\exp(2 \pi i \frac{1}{k})$.
So we get the shift
$$
s^-_c= 2 [\sum_{j=1}^{k-1} (\frac{j}{k}-\frac{1}{2})] =
\frac{k(k-1)}{k}- (k-1) = 0
$$

For an arbitrary $\s$, we regard its cycle decomposition and obtain the result.

\subsection{Theorem} Given a Jacobian Frobenius algebra $A$ up to a twist
by a discrete 
torsion $\alpha \in Z^2(\Sn,k)$ and 
supertwist $\Sigma \in {\rm Hom}(\Sn,\Zz)$ there 
is a unique $\Sn$ Frobenius algebra structure on 
$A^{\otimes n}$.

{\bf Proof.} The uniqueness follows from \S \ref{sn}. The existence result is
deferred to \S \ref{symp} which can be carried over verbatim.

\section{Second quantized Frobenius algebras}
\label{symp}

Given a Frobenius algebra $A$ with multiplication 
$\mu: A\otimes A \rightarrow A$, we can regard its tensor powers 
$T^nA:=A^{\otimes n}$.
These are again Frobenius algebras with the natural tensor multiplication
$\mu^{\otimes n} \in A^{\otimes 3n}=(A^{\otimes n})^{\otimes 3}$, 
tensor metric $\eta^{\otimes n}$ and unit $1^{\otimes n}$.

We can also form the symmetric powers $S^nA$ of $A$. The metric, multiplication
and unit all descend to make $S^nA$ into a Frobenius algebra, but
in terms of general theory [K] we should not regard this object alone,
but rather look at the corresponding orbifold quotient $T^nA/{\mathbb S}_n$.

\subsection{Assumption}
We will assume from now on that $A$ is irreducible and the degree of $A$
is $d$. 

\subsection{Notation} We keep the notation of the previous paragraphs:
$l(\s)$ is the number of cycles in the cycle decomposition of $\s$ and
$|\s|=n-l(\s)$ is the minimal number of transpositions.  

\subsubsection{Lemma} 
Let $\rho$ be the permutation representation of $\Sn$ on $A^{\otimes n}$
permuting the tensor factors. Then the following equations hold
\begin{eqnarray}
\Tr(\rho(\s))&=&\dim(A) l(\sigma)\\
\det(\rho(\s))&=& (-1)^{|\s|{\dim(A)\choose 2}}\begin{cases} 1 & \dim(A) \equiv 0 \text{ or } 1 (4)\\
(-1)^{|\s|} & \dim(A) \equiv 2 \text{ or } 3 (4)
\end{cases}
\end{eqnarray}

{\bf Proof.}
For the first statement we use the fact that 
entries in the standard tensor basis of the matrix
of $\rho(\sigma)$ are just $0$ or $1$. A diagonal entry is $1$ if all of the 
basis elements whose index is in the same subset of $\bar n$ defined by the 
partition $c(\s)$ are equal. The number of such elements is precisely
$\dim(A) l(\s)$.

For the second statement we notice that

$$
\det(\rho(\s))=\det(\rho(\t))^{|\s|}
$$
where $\t$ is any transposition.
For $\tau =(12)$ we  decompose 
$A\otimes A= \bigoplus_{i=1}^{\dim A} e_i \otimes e_i\oplus 
(\bigoplus_{i,j \in \bar n, i \neq j} e_i \otimes e_j) $
for some basis $e_i$ of $A$.
Using this decomposition we find that indeed
$\det(\rho(\s))= (-1)^{|\s|{\dim(A)\choose 2}}$.
For the last statement notice that 
$$
\frac{1}{2}\dim(A)(\dim(A)-1) \equiv  
\begin{cases} 0 (2) &\text{ if } \dim(A) \equiv 0 \text{ or } 1 (4)\\
1 (2) &\text{ if } \dim(A) \equiv 2 \text{ or } 3 (4)
\end{cases}
$$

\subsection{Super-grading} 
As is well known there are only two characters
for $\Sn$: the trivial and the determinant. We will accordingly define 
the {\em parity} with values in $\Z2$
\begin{equation}
\tilde \s \equiv  \begin{cases} 0 \; (2)& \text {if we choose the trivial character}\\
|\s| \; (2)& \text {if we choose the non-trivial character}
\end{cases}
\end{equation}
To unify the notation, we set the parity index $p=0$ in the first case,
which we call even,
and $p=1$ in the second case, which we call odd.

In both cases 
\begin{equation}
\tilde \s = (-1)^{p|\s|}
\end{equation}

\subsection{Intersection algebra structures}

For $\s_1, \dots , \s_m \in \Sn$ we define the following Frobenius algebras:

\begin{eqnarray}
A_{\s}&:=&( A^{l(\s)},\eta^{\otimes l(\s)},1^{\otimes l(\s)})\\
A_{\s_1,  \dots , \s_m}&:=&
( A^{\otimes |\la\s_1,  \dots , \s_m\ra \backslash \bar n|},
\eta^{\otimes | \la\s_1,  \dots , \s_m\ra \backslash \bar n|},
1^{\otimes | \la\s_1,  \dots , \s_m\ra \backslash \bar n|})
\end{eqnarray}

Notice that the multiplication $\mu$ gives rise to a series of maps
by contractions. More precisely given a collection of subsets of $\bar n $
we can contract the tensor components of $A^{\otimes n}$ belonging to the 
subsets by multiplication. Given a permutation we can look at its cycle 
decomposition which yields a decomposition of $\bar n$ into subsets.
We define $\mu(\sigma)$ to be the above contraction. Notice that due to 
the associativity of the multiplication the order in which the contractions 
are
performed is irrelevant.

These contractions have 
several sections. The simplest one being the one mapping 
the product to the first contracted component of each of the 
disjoint contractions.
We denote this map by $j$ or in the 
case of contractions given by $I(\s)$ for some $\s \in \Sn$ 
by $j(\sigma)$. 

E.g.\
$\mu((12)(34))(a\otimes b\otimes c \otimes d)= ab\otimes cd$ 
and $j((13)(24))(ab\otimes cd) = ab \otimes cd \otimes 1\otimes 1$.

Thus we define the following maps
\begin{eqnarray}
r_{\s}:A_e \rightarrow A_{\s}&;& r_{\s}:=\mu(\sigma)\\
i_{\s}: A_{\s}\rightarrow A_e&;& i_{\s}:=j(\sigma)
\end{eqnarray}

Moreover the same logic applies to the spaces $A_{\s_1, \dots, \s_m}$
and we similarly define
$r_{\s_1, \dots, \s_m}, i_{\s_1, \dots, \s_m}$
where the indices are symmetric
and maps
\begin{equation}
r_{\s_1, \dots, \s_m}^{\s_1, \dots, \s_{m-1}}:
A_{\s_1, \dots, \s_{m-1}}\rightarrow 
A_{\s_1, \dots, \s_m} , 
i_{\s_1, \dots, \s_m}^{\s_1, \dots, \s_{m-1}}:
A_{\s_1, \dots, \s_m}\rightarrow A_{\s_1, \dots, \s_{m-1}}
\end{equation}

where the again the indices are symmetric.

We also notice that $A_{\s} = A_{\s^{-1}}$ and $A_{\s,\s}=A_{\s}$.

\subsection{Remark}
The sections $i_{\s}$ also satisfy the condition 
\begin{equation}
i_{\s}(ab_{\s})=\pi(a)i_{\s}(b_{\s})
\end{equation}

\subsection{Proposition} 
\label{traces}
The maps $r_{\sigma}$ make
$A_{\s},\eta_{\s}$ into a special $\Sn$ reconstruction data.
A  choice of parity $\tilde \s$ fixes the character to be: 
\begin{equation}
\chi_{\s}=(-1)^{p|\s|}
\end{equation}

Furthermore the collection of maps 
$r_{\s_1, \dots, \s_m}^{\s_1, \dots, \s_{m-1}}$ turns the collection
of $A_{\s_1, \dots, \s_m}$ into special intersection $\Sn$ reconstruction data.

{\bf Proof.} It is clear that all the $A_{\s}$ are cyclic $A_e$ modules
 and is is clear
that $A_{\sigma} = A_{\sigma^{-1}}$.

Also the $\eta_{\s}$ remain unscaled since
$(-1)^{p|\s|}\chi_{\s}\equiv 1$.

What remains to be shown is that the character 
is indeed given by $\chi_{\s}=(-1)^{p|\s|}$ and that the trace axiom holds.

This is a nice exercise. We are in the graded case and moreover the identity
is up 
to scalars the only  element with degree zero ---
unless ($\dim A=1$) and we are in the
case of $pt/\Sn$ which was considered in \ref{point}. 
So if $c\in A_{[\s,\s']}: c \neq \lambda 1_e$
then the trace axiom is satisfied automatically.

Therefore we only need to consider the case $c= 1 \in A_{[\s,\s']}$ with 
$[\s,\s']=e$.
In this case, we see that $\s'$ acts on $A_{\s}\iso A^{\otimes l(\s)}$ 
as a permutation. Indeed the normalizer of $\s$ 
is the semi--direct product of permutations of the cycles and cyclic 
groups whose induced action on
$A_{\s}$ is given by permutation and identity respectively. 

We claim the trace has the value 
\begin{equation}
\Tr \varphi_{\s}|_{A_{\s'}}=\dim(A_{\s,\s'})
\end{equation}

This is seen as follows. Looking at the permutation action on the factors 
of $A_{\s}$, we see that the trace has entries 0 and 1 in any fixed basis of 
$A_{\s}$ induced by a fixed choice of basis of $A$. The value 1 appears
if the pure tensor element has exactly the same entry in all tensor 
components labelled by elements which are in the same cycle of $\s$ (acting on 
$A_{\s}$). But these are precisely the elements that span $A_{\s,\s'}$. 
To be more precise there is a canonical isomorphism of these element with
$A_{\s,\s'}$ given by tensors of iterated diagonal maps $\Delta: A \rightarrow
A\otimes A, \Delta(a) = a\otimes a$.

Thus the trace axiom can be rewritten as:
\begin{equation}
\label{tracesigma}
\chi_{\s}\varphi_{\s,\s'}(-1)^{p|\s|}= \chi^{-1}_{\s'}
\varphi_{\s',\s}(-1)^{p|\s'|}
\end{equation}

In particular if $\s'=e$
$$
(-1)^{p|\s|} \dim(A_{\s})= \chi_{\s} \Tr(\rho(\s)|_{A^{\otimes n}})
$$
so that 
$$
\chi(\s)= (-1)^{p|\s|}
$$

Combining the above we find that:
\begin{equation}
\label{tracevalue}
\chi_{\s} \STr(\phi_{\s}|_{A_{\s'}})= 
(-1)^{p(|\s||\s'|+|\s|+|\s'|)}\dim(A_{\s,\s'})
\end{equation}
which is an expression completely symmetric in $\s,\s'$ and invariant
under a change $\s \mapsto \s^{-1}$.

For the last statement we only need to notice that consecutive contractions
yield commutative diagrams which are co--Cartesian.
The structural isomorphisms being clear since they
can all be given by the identity morphism --- there is no rescaling.

\subsection{Proposition (Algebraic Discrete Torsion)}
Fix the $\sign \equiv 1$ and $\sign^{\s}\equiv 1$ and set
$(-1)^{\tilde \s^{\s'}} = \det_{V{\s'}}(\s)
=(-1)^{\codim_{V_{\s'}}(V_{\s,\s'})}$ 
where $ \det_{V{\s'}}(\s)$ 
is the determinant of the induced action of $\s$ on the fixed point
set of $\s'$. Furthermore fix $\chi_{\s}^{\s'}$ by 
$(-1)^{p(\codim_{V_{\s'}}(V_{\s,\s'}))}$. 
Then $\sign$ and the $\sign^{\s}$ are compatible  
and 

\begin{equation}
\eps(\s,\s')= (-1)^{p(|\s||\s'|)}(-1)^{p|\s|}
(-1)^{p(\codim_{V_{\s'}}(V_{\s,\s'}))}
\end{equation}
or in the notation of \ref{cent}

$$
\eps(\s',\s)=\begin{cases} (-1)^{p(k|\s|+k+1)}&
\text{if $\s'=\t_k$} \\
(-1)^{p((k-1)|\s|+(k-1))}& 
\text{if $\s'=c_k$}
\end{cases}
$$

{\bf Proof.} First:
$$
\nu(\s,\s') \equiv \codim_{V_{\s'}}(V_{\s,\s'}) +\codim(V_{\s}) (2)
$$ 
which satisfies \ref{nucompat}, since
$$
\codim_{V_{\s}}(V_{\s,\s'}) +\codim(V_{\s})= \codim(V_{\s,\s'})
=\codim_{V_{\s'}}(V_{\s,\s'}) +\codim(V_{\s'})
$$
Now just by definition
$$
\eps(\s,\s')= (-1)^{p(|\s||\s'|)}(-1)^{p|\s|}
(-1)^{p(\codim_{V_{\s'}}(V_{\s,\s'}))}
$$ 
and lastly: $\codim_{V_{\s}}(V_{\t_k,\s})=1$ and 
$\codim_{V_{\s}}(V_{c_k,\s})=0$.

\subsubsection{Remark} This algebraic discrete torsion indeed reproduces
the effect that turning it on yields the super--structure on the
twisted sectors as postulated in [D2]. The computation of the discrete
torsion in [D2] was however done for $pt/\Sn$ with the choice of
cocycle $\g$ given by a Schur multiplier, see \ref{schurdt}.
The current calculation explains how the non--trivial Schur--multiplier
used to twist by a discrete torsion behaves like a supertwist.
In terms of  \ref{qdef} one can see this as the
fact that in both twists --super and non--trivial discrete torsion--
$q=1$.

\subsection{Proposition}
\label{uniqueness} 
After possibly twisting by discrete torsion 
any  cocycle $\g$ compatible with the
special reconstruction data is normalizable and hence unique after the 
normalization.

{\bf Proof.} Verbatim the proof of \ref{normalizable}.

So from now on we can and will deal with normalized cocycles.

\subsubsection{Lemma}
\label{pullback}
For any  minimal decomposition $T$ of $\s'$ into transpositions 
$\s'= \t_1 \dots \t_{|\s'|}$
\begin{equation}
\check r_{\s}(1_{\s})= \prod_{i} \g_{\t_i,\t_i} 
\end{equation}

{\bf Proof.} Notice that $I_{\s}= \bigoplus I_{\t_i}$
and thus $I_{\s} \prod_{i \in I} \g_{\t_i} =0$.
Furthermore \\
$\deg(\prod_{i \in I} \g_{\t_i})=d|\s|= s^+(\s)
= 2d_{\s} = \deg(\g_{\s,\s^{-1}})$
and 
$\dim(I_{\d})^{dl(\s)}=\dim(A^{\otimes n})-1$ 
where the  superscript  denotes the 
part of homogeneous degree. This follows from the equalities:
$\dim((I_{\s})^{dl(\s)}) = \dim (\mathrm{Ker} (r_{\s}^{dl(\s)})) 
= \dim(A^{\otimes n})-\dim(\mathrm{Im}( r_{\s}^{dl(\s)}) )
=\dim(A^{\otimes n}) -1$.
We split $(A^{\otimes n})^{dl(\s)} = (I_{\s})^{dl(\s)} \oplus L$
where $L$ is the line generated by $i_{\s}(\rho_{\s})$.

We have to show that
$$
\eta(\prod_{i \in I} \g_\t,b) = \eta_{\s}(1_{\s}r_{\s}(b))
$$

This is certainly true if $\deg(b) \neq dn-d|\s|=dl(\s)$ since then both 
sides vanish. This is also the case if $b\in I_{\s}$. It remains to show 
that
$\eta(\prod_{i \in I} \g_{\t_i,\t_i},i_{\s}(\rho_{\s}))=1$. 

We do this by induction on $|\s|$, the statement being clear for $|\s|=1$.
Let $\t_{|\s|}= (ij)$ and set $\s'=\s\t_{|\s|}$ then 
$$
i_{\s'}(\rho_{\s'}) = \g_{\t_{|\s|},\t_{|\s|}}i_{\s}(\rho_{\s})
$$
which follows from the equation
$\rho\otimes 1 \g_{(12),(12)} = \rho \otimes \rho$ and its pull back.
So
$$1=\eta(\prod_{i=1}^{|\s|-1} \g_{\t_i,\t_i}, i_{\s'}(\rho_{\s'}))
= \eta(\prod_{i=1}^{|\s|-1} \g_{\t_i,\t_i},\g_{\t_{|\s|},\t_{|\s|} }
i_{\s}\rho_{\s})
$$
Another way to see this is to use 
the isomorphism $A_{\s} \iso A_{\t_1,\dots,\t_{|\s|}}$ and the iterated
restriction maps for the pull-back, noticing, that indeed the 
$\g_{\t,\t}$ pull back onto each other in the various space.

Using the same rationale we obtain:

\subsection{Corollary}
\label{intersectionpullback}
\begin{equation}
 \check
r_{\s,\s'}^{\s\s'} (1_{\s,\s'})=
\pi_{\s\s'}(\prod_{i\in I_{\s,\s'}}
\g_{\t_i,\t_i})
\end{equation}
where $I_{\s,\s'}=
\{i\in I: |\la \s\s',\t \ra \backslash \bar n|< 
|\la \s\s' \ra \backslash \bar n|\}$ or in other words the $\g_{\t_i,\t_i}$ 
that do not get contracted.

\subsection{Grading and shifts}
The meta--structure for symmetric powers is given by treating $A^n$ as the
linear structure, just like the variables in the Jacobian case.
In particular we fix the following degrees and shifts
\begin{equation*}
\deg(1_{\s})=d|\s| 
\end{equation*}
\begin{eqnarray*}
s^+_{\s}&=& d|c(\s)| , \quad s^-_{\s}= 0 \\
s_{\s} &=&\frac {1}{2}(s_{\s}^+ + s_{\s}^-)= \frac{d}{2}|c(\s)|
\end{eqnarray*}

Notice that as always there is no ambiguity for $s^+$, not even in the 
choice of dimension of $A_{\s}$, but the choice for $s^-$ is 
a real one which is however the only choice which extends the natural grading
if $A$ is Jacobian.

This view coincides with the realization of $A^{\otimes n}$ as
the n-th tensor product of the extension of 
coefficients to A of the Jacobian algebras for $f=z^2$.

\subsection{Notation}
The geometry of $\Sn$--Frobenius algebras is given by the subspace arrangement
of fixed point sets $V_{\s}= Fix(\s)\subset k^n$ of the various $\s \in \Sn$ acting on $k^n$ as
well as their intersections $V_{\s,\s'}= V_{\s}\cap V_{\s'}$, etc.,
which were introduced in  \S \ref{sn}. 

Recall that $|\s|= \mathrm{codim}_V(V_{\s})$. We also define
$|\s,\s'|:=  \mathrm{codim}(V_{\s}\cap V_{\s'})$
and set

\begin{eqnarray}
d_{\s,\s'}&:=& \frac{1}{d}\deg(\g_{\s,\s'})=
\frac{1}{2} (|\s|+|\s'|-|\s\s'|)\nn\\
n_{\s,\s'}&:=&\frac{1}{d}\deg(\check {r}_{\s,\s'}^{\s\s'}(1_{\s,\s'}))
=\codim_{V_{\s\s'}}(V_{\s,\s'})
\nn\\
&=&|\s,\s'|-|\s\s'|\nn\\
\tilde g_{\s,\s'}&:=&\frac{1}{d}\deg(\tilde {\g}_{\s,\s'})
= d_{\s,\s'}-n_{\s,\s'}\nn\\
&=&\frac{1}{2}(|\s|+|\s'|+|\s\s'|-2|\s,\s'|)
\end{eqnarray}

Now given two elements $\s,\s'\in \Sn$ their representation on
$k^n$ naturally splits $k^n$ into a direct sum, which is given by the smallest 
common block decomposition of both $\s$ and $\s'$. More precisely:

Fix the standard basis $e_i$ of $k^n$. 
For a subset $B \in \bar n$ we set 
$V_{B}= \bigoplus_{i\in B} k e_i \subset k^{n}$.
Given $\s,\s'$ we decompose
$$
V:= k^{n }= \bigoplus_{B\in \la\s,\s' \ra \backslash \bar n}  V_{B}
$$

and decompose
\begin{equation}
V_{\s} = \bigoplus_{B\in \la\s,\s' \ra \backslash \bar n}  V_{\s;B}; \quad
V_{\s,\s} = \bigoplus_{B\in \la\s,\s' \ra \backslash \bar n}  V_{\s,\s;B}
\end{equation}
where $V_{\s;B}:= V_{\s}\cap V_B; V_{\s,\s';B}:= V_{\s,\s'}\cap V_B $ 
and we used the notation of \ref{snnotation}.

Notice that $\dim(V_{\s,\s';B})=1$ and we can decompose $\tilde\g_{g,h}=\bigotimes_B \tilde \g_{g,h;B}$.

Using the notation:
$$
|\s|_B := \codim_{V_B} (V_{\s;B}), \quad 
|\s,\s'|_B := \codim_{V_B} (V_{\s,\s;B})
$$
set
\begin{eqnarray}
d_{\s,\s';B}&:=&\frac{1}{2}(|\s|_B+|\s'|_B-|\s,\s'|_B)\nn\\
n_{\s,\s';B}&:=&|\s,\s'|_B-|\s\s'|_B= \codim_{V_{\s\s';B}}(V_{\s,\s';B})\nn\\
\tilde g_{\s,\s';B}&:=& d_{\s,\s';B}- 
n_{\s,\s';B}=\frac{1}{d}\deg(\tilde {\g}_{\s,\s';B})\nn\\
&=&\frac{1}{2}(|\s|_B+|\s'|_B+|\s\s'|_B-2|\s,\s'|_B)
\end{eqnarray}
Notice that all the above functions take values in ${\bf N}$.

\subsubsection{Triple intersections}
For any number of elements $\s_i$ we can analogously define the 
above quantities. We will do this for the triple intersections,
since we need these to show associativity and although tedious we 
do this in order to fix the notation.

We regard the triple intersections $V_{\s,\s',\s''}= V_{\s}\cap V_{\s'}\cap
V_{\s''}$.

Recall that $|\s|= \mathrm{codim}_V(V_{\s})$. We also define
$|\s,\s',\s''|:=  \mathrm{codim}(V_{\s,\s',\s''})$
and set

\begin{eqnarray}
d_{\s,\s',\s''}&:=& \frac{1}{d}\deg(\g_{\s,\s'}\g_{\s\s',\s''})\nn \\
&=&\frac{1}{2} (|\s|+|\s'|-|\s\s'|+|\s\s'|+|\s''|-|\s\s'\s''|)\nn\\
&=&\frac{1}{2} (|\s|+|\s'|+|\s''|-|\s\s'\s''|)\nn\\
n_{\s,\s',\s''}&:=&
\frac{1}{d}\deg(\check {r}_{\s,\s',\s''}^{\s\s'\s''}(1_{\s,\s'\s''}))
= \codim_{V_{\s\s'\s''}}(V_{\s,\s',\s''})\nn\\
&=&|\s,\s',\s''|-|\s\s'\s''|\nn\\
\tilde g_{\s,\s',\s'}&:=&\frac{1}{d}\deg(\tilde {\g}_{\s,\s',\s''})= d_{\s,\s',\s''}-n_{\s,\s',\s''}\nn\\
&=&\frac{1}{2}(|\s|+|\s'|+|\s''|+|\s\s'\s''|-2|\s,\s',\s''|)
\end{eqnarray}

where $\tilde \g_{\s,\s',\s''}$ was defined in (\ref{tildetripel}).

As above given three elements $\s,\s',\s''\in \Sn$ their representation on
$k^n$ naturally splits $k^n$ into a direct sum, which is given by the smallest 
common block decomposition of $\s,\s'$ and $\s'$. More precisely:

Again, fix the standard basis $e_i$ of $k^n$. 
For a subset $B \in \bar n$ we set 
$V_{B}= \bigoplus_{i\in B} k e_i \subset k^{n}$.
Given $\s,\s'$ we decompose
$$
V:= k^{n }= \bigoplus_{B\in \la\s,\s',\s'' \ra \backslash \bar n}  V_{B}
$$
and decompose
\begin{eqnarray}
&V_{\s}& = \bigoplus_{B\in \la\s,\s',\s'' \ra \backslash \bar n}  V_{\s;B}; 
\;
V_{\s,\s'} 
= \bigoplus_{B\in \la\s,\s',\s'' \ra \backslash \bar n} V_{\s,\s';B}\nn\\
&V_{\s,\s',\s''}&
=\bigoplus_{B\in \la\s,\s',\s'' \ra \backslash \bar n} V_{\s,\s',\s'';B}
\end{eqnarray}
where $V_{\s;B}:= V_{\s}\cap V_B; V_{\s,\s';B}:= V_{\s,\s'}\cap V_B
V_{\s,\s',\s'';B}:= V_{\s,\s',\s''}\cap V_{B} $

Notice that $\dim(V_{\s,\s',\s'';B})=1$. 

We will also use the notation:
$$
|\s|_B := \codim_{V_B} (V_{\s;B}), 
|\s,\s'|_B := \codim_{V_B} (V_{\s,\s;B})$$ 
and
$$
|\s,\s',\s''|_B := \codim_{V_B} (V_{\s,\s,\s'';B})
$$

\subsection{The cocycle in terms of $\g_{\t,\t}$s}
\label{cocycle}

Let $\g_{\s,\s'}$ be given by the following:
 
For transversal $\s,\s'$ we set $\g_{\s,\s'}=1$.

If
$\s$ and
$\s'$ are {\em not} transversal 
using Theorem \ref{normalize} we set 

\begin{eqnarray}
\label{nottrans}
r_{\s,\s'}(\g_{\s,\s'})&=&r_{\s\s'}(\prod_{i\in I}\g_{\t_{i},\t_{i}})=
\prod_{i\in I'}\pi_{\s\s'}(\g_{\t_{i},\t_{i}})\prod_{j\in I''}
r_{\s\s'}(\g_{\t_{j},\t_{j}})\nn\\
&=:&\bar\g_{\s,\s'}\g_{\s,\s'}^{\perp}
\end{eqnarray}
where
\begin{eqnarray}
I'=\{ i \in I:
\pi_{\s,\s'}(\g_{\t_{i},\t_{i}})=\pi_{\s\s'}(\g_{\t_{i},\t_{i}})\}\nn\\
I''=\{ i \in I: \pi_{\s,\s'}(\g_{\t_{i},\t_{i}})
\neq\pi_{\s\s'}(\g_{\t_{i},\t_{i}})\}
\end{eqnarray} 
and
$\bar\g_{\s,\s'}\in i_{g,h}^{gh}(A_{g,h})$

\subsection{Proposition}
\label{existence}
The equations of \ref{cocycle} are well defined and yield a group cocycle
compatible with the reconstruction data.
Furthermore 
\begin{eqnarray}
\g_{\s,\s'}^\perp&=& r_{\s,\s'}^{\s\s'}(1_{\s,\s'})\label{1}\\
\bar\g_{\s,\s'}&=&i_{\s,\s'}^{\s\s'}
(\bigotimes_{B\in \la \s,\s'\ra \backslash \bar n} e^{g(\s,\s',B)})\label{2}\\
\g_{\s,\s'}&=& r_{\s,\s'}^{\s\s'}
(\bigotimes_{\in \la \s,\s'\ra \backslash \bar n} e^{g(\s,\s',B)})=
\check r_{\s,\s'}^{\s,\s'}(\tilde \g_{g,h})
\label{3}
\end{eqnarray}

{\bf Proof.}
We need to check that indeed equation (\ref{nottrans}) is well defined.
From Lemma \ref{pullback} and the Corollary \ref{intersectionpullback}
we know that (\ref{1}) 
is true and that the  product over $I''$ is well defined.

For (\ref{2}) we notice that if a $\g_{\t_i,\t_i}$ gets contracted,
then 
\begin{equation}
\label{Euler}
\pi_{\s\s'}(\g_{\t_i,\t_i})= 1 \otimes \dots \otimes 1\otimes e \otimes 1
\otimes \dots \otimes 1
\end{equation}
where 
$e = \mu \check\mu(1)$ is the Euler class which sits in the image of the  
$k$-th factor which is the same as the image of the $l$-th factor 
under the map $\pi_{\s,\s'}$ if $\t_i =(kl)$.
%Here we used $\Delta := \check\mu: A\rightarrow A\otimes A$.

The well definedness then follows by decomposition 
into $V_{\s\s',B}$ from the statement for one--dimensional 
$V_{\s,\s'}$  where it is clear from grading.

Finally (\ref{3}) follows from (\ref{1}) and (\ref{2}) via 
Proposition \ref{intersect}.

For the associativity we use the general theory of intersection
algebras \ref{ass}.
Here we notice that indeed the number of 
$\g_{\t_i,\t_i}:i \in I''$ contracted in each component $B$
by $r_{\s\s',\s''}^{\s,\s'}$ is given by 
\begin{multline}
 n_{\s,\s';B}- |\s,\s',\s''|_B +|\s\s',\s''|_B\\
= |\s,\s'|_B-|\s\s'|_B- |\s,\s',\s''|_B +|\s\s',\s''|_B 
:=q(\s,\s',\s'';B)
\end{multline} 
so that by commutativity of (\ref{assdiagram})
$$
r_{\s\s',\s''}^{\s\s'}(\g_{\s,\s'})= 
i_{\s,\s',\s''}^{\s\s',\s''}(\bigotimes_B e^{q(\s,\s',\s'';B)})
\check r_{\s,\s',\s''}^{\s\s',\s''}(1_{\s,\s',\s''})
$$
and thus the $\g^{\perp}$ match also by commutativity.
What remains to be calculated is the power of $e$ in each of the
components $B$.
This power is given by
\begin{eqnarray*}
&&\frac{1}{2}(|\s\s'|_B+|\s''|_B+|\s\s'\s''|_B-2|\s\s',\s''|_B)\\
&+&\frac{1}{2}(|\s|_B+|\s'|_B+|\s\s'|_B-2|\s,\s'|_B)\\
&+&(|\s,\s'|_B-|\s\s'|_B- |\s,\s',\s''|_B +|\s\s',\s''|_B )\\
&=&\frac{1}{2}(|\s|_B+|\s'|_B+|\s''|_B+|\s\s'\s''|_B-2|\s,\s',\s''|_B)\\
&=&\tilde g_{\s,\s',\s'';B}
\end{eqnarray*}
q.e.d.
\medskip

Putting together the Propositions \ref{uniqueness} and
\ref{existence} of this section we obtain:

\subsection{Theorem} There exists a unique normalized cocycle compatible
with the above special reconstruction data. There is only one compatible
cocycle in the all even case.

In the super--case there are
two choices of parity for the twisted sectors: all even or the parity of
$p|\s| \equiv |\s| (2)$.  Fixing the parity fixes the non--abelian cocycle.

In other words, there is a unique multiplicative $\Sn$ Frobenius
algebra structure on the tensor powers of $A$ 
and there are two $G$--actions labelled by parity.

\subsection{Definition} We call the {\it symmetric power of a
Frobenius algebra}
 the $\Sn$--twisted Frobenius algebra obtained from $T^nA, (r_{\s})$ by
using the unique normalized cocycle with all even sectors
and the {\it super--symmetric power of a
Frobenius algebra}
 the $\Sn$--twisted Frobenius algebra obtained from $T^nA, (r_{\s})$ by
using the unique normalized cocycle with the  parity given by
$A_{\s} \equiv |\s| \; (2)$.

\subsection{Definition} We  define {\em the 
second quantization of a Frobenius algebra $A$} to be the sum of all
symmetric powers of $A$ and
{\em the 
second super--symmetric quantization of a Frobenius algebra $A$} 
to be the sum of all
super--symmetric powers of $A$. 
We consider this sum either as formal or as a direct sum,
where we need to keep in mind that the degrees of the summands are not equal.

\subsection{Comparison with the Lehn and Sorger construction}

In [LS] Lehn and Sorger constructed a non--commutative multiplicative
structure  in the special setting of symmetric powers. 
By the uniqueness result of the last
section we know ---since their cocycles are also normalized--- that
their construction has to agree with ours. In this section we make this
explicit. Our general considerations of intersection algebras
explain the appearance of their
cocycles as the product over the Euler class to the graph defect times
contribution stemming from the dual of the contractions.

\subsubsection{Definition}(The graph defect) For
$ B\in \la\s,\s' \ra \backslash \bar n$  define the 
graph defect as [LS]
 
\begin{equation}
g(\s,\s';B):= \frac{1}{2}
(|B|+2 - 
|\la\s\ra \backslash B|- |\la\s'\ra \backslash B|-
|\la\s,\s'\ra \backslash B|)
\end{equation}

The equality of the two multiplications follows from:

\subsubsection{Proposition} 
$$
g(\s,\s';B)=\tilde g_{\s,\s';B}=
\frac{1}{2}(|\s|_B+|\s'|_B+|\s\s'|_B-2|\s',\s'|_B)
=d_{\s,\s';B}-n_{\s,\s';B}
$$

{\bf Proof.} 
By the above:
\begin{eqnarray*}
g(\s,\s';B)&=& \frac{1}{2}(
\dim(V_B)+2\dim_{V_B}(V_{\s,\s';B})- \dim(V_{\s,B})-\\
&&-\dim(V_{\s';B})-\dim(V_{\s\s';B}))\\
&=&\frac{1}{2}(\dim(V_B)- \dim(V_{\s,B}) + \dim(V_B)-\dim(V_{\s';B})
\\
&&+\dim(V_B) -\dim(V_{\s\s';B})-(2\dim(V_B)- \dim_{V_B}(V_{\s,\s';B})))\\
&=&\frac{1}{2}(|\s|_B +|\s'|_B +|\s\s'|_B-2|\s,\s'|_B)
\end{eqnarray*}

\subsubsection{Remark} The above equation makes it obvious that $g \in \mathbb{N}$, since $d_{\s,\s';B},$ $n_{\s,\s';B}\in \mathbb{N}$ and both
$|\s|+|\s'|\geq |\s,\s'|$ and $|\s\s'|\geq |\s,\s'|$. The first inequality follows from $V_{\s,\s'}=V_{\s}\cap V_{\s'}$ and the second one from 
$V_{\s,\s'}\subset V_{\s\s'}$. 

\subsection{Remark} The change of sign needed to recover the cohomology algebra
of the Hilbert scheme of a K3 surface can also be obtained by a twisting
with a discrete torsion. To be precise by the normalized discrete torsion class $\a$ defined by $\a(\t,\t)=-1$ ($\t\in \Sn,|\t|=1$), see \cite{K5}.

\renewcommand{\theequation}{A-\arabic{equation}}
\renewcommand{\thesection}{A}
% redefine the command that creates the equation no.
\setcounter{equation}{0}  % reset counter 
\setcounter{subsection}{0}
\section*{Appendix A}  % use *-form to suppress numbering

\subsection{Theorem}
\label{normal}
Any normalizable  non--abelian $\Sn$ cocycle 
$\varphi$ with values in $k^{*}$ can be normalized after a 
rescaling and then one of the following holds:
 $\forall \s, \t, |\t|=1: $
 $$
 \varphi_{\s,\t}= 1
 $$
We call this case even and set the parity $p=0$.
  Or 
 $\forall \s,\t |\t|=1$
 $$
 \varphi_{\s,\t}= (-1)^{|\s|}
 $$
We call this case odd and set the parity $p=1$.
In unified notation:
\begin{equation}
\varphi_{\s,\t}= (-1)^{p|\s|}
\end{equation}
with $p\in \{0,1\}$.  

{\bf Proof.}  
By assumption
\begin{equation}
\label{commuting}
\forall \t,\t',|\t|=|\t'|=1,[\t,\t']= e: \varphi_{\t,\t'}=(-1)^p
\end{equation}

We will show by induction that we can scale such that
\begin{equation}
\label{notcommuting}
\forall \t,\t',|\t|=|\t'|=1,[\t,\t']\neq e: \varphi_{\t,\t'}=(-1)^p
\end{equation}

Combining (\ref{commuting}) and (\ref{notcommuting}):
\begin{equation}
\label{transpos}
\forall \t,\t';|\t|=|\t'|=1: \varphi_{\t,\t'}=(-1)^p
\end{equation}

{\sc Induction for (\ref{transpos}).}

Assume that (\ref{transpos}) holds for $\t,\t' \in \Sn \subset \Snn$.

Now scale with 
\begin{eqnarray}
\l_{(ij)}&:=&(-1)^p \varphi_{(n-1 \, n+1),(n\, n+1)} 
\text{ for } i,j \leq n\nn\\ 
\l_{(i \, n+1)}&:=& (-1)^p \varphi_{(in),(n \, n+1)} \text{ for } i <n
\nn\\
\l_{(n \, n+1),(n\, n+1)}&:=&1
\end{eqnarray}

Notice this implies that
\begin{eqnarray}
\label{n+1}
\tilde \varphi_{(in),(n \, n+1)}
&=&\frac{\l_{(n\, n+1)}}{\l_{(i\, n+1)}} \varphi_{(in),(n\, n+1)}\nn\\
&=&(-1)^p \frac{1}{\varphi_{(in),(n \,n+1)}} \varphi_{(in),(n\, n+1)}
=(-1)^p\\
\tilde \varphi_{(in),(i\, n+1)}&=&\tilde \varphi^{-1}_{(in),(n\, n+1)}=(-1)^p\\
\label{in}
\tilde \varphi_{(ij),(kl)}&=&\varphi_{(ij),(kl)}=(-1)^p 
\text { if } i,j,k,l \leq n
\end{eqnarray}
where the last statement follows by induction.

We need to show
\begin{equation}
\tilde \varphi_{\t,\t'}=(-1)^p
\end{equation}

For $n=2$ the statement is true.

So we assume $n\geq 2$ and by assumption:
\begin{equation}
\label{comsol}
\forall \t,\t';|\t|=|\t'|=1;[\t,\t']=e:\quad \varphi_{\t,\t'}=(-1)^p
\end{equation}

Thus by induction (\ref{n+1})--(\ref{in}) and (\ref{comsol}),
we need to check the cases

\begin{itemize}
\item[i)] $\t=(ij),\t'=(j\, n+1); i,j \in \{1,\dots,n-1\}; i\neq j$
\item[ii)] $\t=(i\, n+1),\t'=(j\,n+1); i,j \in \{1,\dots,n\}, i\neq j$
\item[iii)] $\t=(in+1),\t'= (ij);i,j\in \{1,\dots,n\}, i\neq j$
\end{itemize}

Notice that 
$\tilde \varphi_{(i\, n+1),(ij)}=\tilde \varphi_{(i\, n+1),(j\, n+1)}^{-1}$ 
and thus ii) implies iii). Else
iii) follows by (\ref{n+1}) and thus it suffices to show i) and ii).

For i)
\begin{eqnarray*}
\tilde \varphi_{(ij),(j\, n+1)}&=&\frac{\l_{(j\, n+1)}}{\l_{(i\, n+1)}}
\varphi_{(ij),(j\, n+1)}
=\frac{\varphi_{(jn)(nn+1)}}{\varphi_{(in)(nn+1)}} \varphi_{(ij),(j\, n+1)}\\
&=&\varphi_{(jn),(n\, n+1)}\varphi_{(in),(i\, n+1)} \varphi_{(ij),(j\, n+1)}
=\varphi_{(in)(ij)(jn),(n\, n+1)}\\
&=&\varphi_{(ij),(n\, n+1)}=(-1)^p
\end{eqnarray*}
by (\ref{comsol}).

for ii)
If $j=n$ then
\begin{eqnarray*}
\tilde \varphi_{(i\, n+1),(n\, n+1)}&=&\frac{\l_{(n\, n+1)}}{\l_{(ij)}}
\varphi_{(i\, n+1),(n\, n+1)}
=\frac{(-1)^p}{\varphi_{(n-1\, n+1)(n\, n+1)}} \varphi_{(i\, n+1),(n\, n+1)}
\end{eqnarray*}
so if $i=n-1$ $\tilde \varphi_{(n-1\, n+1),(n\, n+1)}=(-1)^p$

If $i\neq n-1$ then
\begin{eqnarray*}
(-1)^p\frac{\varphi_{(i\, n+1),(n\, n+1)}}{\varphi_{(n-1\, n+1)(n\, n+1)}} 
&=&(-1)^p\varphi_{(n-1\, n+1),(n-1\, n)}\varphi_{(i\, n+1),(n\, n+1)}\\
&=&(-1)^p\varphi_{(i\, n+1)(n-1\,n+1),(n-1\, n)}\\
&=&(-1)^p\varphi_{(n-1\, i)(i\,n+1),(n-1\, n)}\\
&=&(-1)^p\varphi_{(n-1\, i),(n-1\,n)}\varphi_{(i\,n+1),(n-1\, n)}=(-1)^p
\end{eqnarray*}

If $j\neq n$
\begin{eqnarray*}
\tilde \varphi_{(i\, n+1),(j\, n+1)}&=&\frac{\l_{(j\, n+1)}}{\l_{(ij)}}
\varphi_{(i\, n+1),(j\, n+1)}
=\frac{\varphi_{(jn)(n\, n+1)}}{\varphi_{(n-1\, n+1)(n\, n+1)}} 
\varphi_{(i\, n+1),(j\, n+1)}\\
&=&\varphi_{(jn),(n\, n+1)} \varphi_{(n-1\, n+1),(n-1\, n)} 
\varphi_{(i\, n+1),(j\, n+1)}\\
&=&\varphi_{(i\, n+1)(jn)(n-1\, n+1),(n-1\, n)}
\end{eqnarray*}

Now first assume $\{i,j\}\cap \{n-1,n\}=\emptyset$ then

\begin{eqnarray*}
\varphi_{(i\, n+1)(jn)(n-1\, n+1),(n-1\, n)}
&=&\varphi_{(n+1\, n-1) (n-1\, i)(jn),(n-1\, n)}\\
&=&\varphi_{(jn),(n-1n)}\varphi_{(n-1\, i),(n-1\, j)}
\varphi_{(n-1\, n+1),(ij)}
=(-1)^p\\
\end{eqnarray*}

Case 2a) $i=n,j=n-1$ then 
\begin{eqnarray*}
\varphi_{(i\, n+1)(jn)(n-1\, n+1),(n-1\, n)}
&=&\varphi_{(n\, n+1)(n-1n)(n-1\, n+1),(n-1\, n)}\\
&=&\varphi_{(n-1\, n),(n-1\, n)}=(-1)^p
\end{eqnarray*}

Case 2b) $i=n, j\neq n-1$
\begin{eqnarray*}
\varphi_{(i\, n+1)(jn)(n-1\, n+1),(n-1\, n)}
&=&\varphi_{(n\, n+1)(jn)(n-1\, n+1),(n-1\, n)}\\
&=&\varphi_{(n-1\, n)(jn)(j\, n+1),(n-1\, n)}\\
&=&\varphi_{(n-1\, n),(j\,n-1)} \varphi_{(jn), (n-1\, n)}
\varphi_{(j\, n+1),(n-1\, n)}=(-1)^p
\end{eqnarray*}

Case 3) $i=n-1, j\neq n$
\begin{equation*}
\varphi_{(i\, n+1)(jn)(n-1\, n+1),(n-1\, n)}
=\varphi_{(n-1\, n+1)(jn)(n-1\, n+1),(n-1\, n)}=\varphi_{(jn),(n-1\,n)}=(-1)^p
\end{equation*}

Case 4) $j=n-1, i\neq n$

\begin{eqnarray*}
\varphi_{(i\, n+1)(n-1\, n)(n-1\, n+1),(n-1\, n)}
&=&\varphi_{(n\, n-1)(n-1\, i)(i\, n+1)),(n-1\, n)}\\
&=&\varphi_{(n\, n-1)(in)}\varphi_{(n-1\, i),(n-1\, n)}\varphi_{(i\, n+1)),(n-1\, n)}=(-1)^p
\end{eqnarray*}

Finally if $\s=\prod_{i=1}^{|\s|}\t_i$
$$
\varphi_{\s,\t}= \prod_{i=1}^{|\s|}\varphi_{\t_i,\tilde \t_i}= (-1)^{p|\s|}
$$
with $\tilde \t_i= (\prod_{j=i+1}^{|\s|}\t_i) 
\t(\prod_{j=i+1}^{|\s|}\t_i)^{-1}$,
$|\tilde \t_i|=|\t|=1$
which proves the Theorem.

\renewcommand{\theequation}{B-\arabic{equation}}
\renewcommand{\thesection}{B}
% redefine the command that creates the equation no.
\setcounter{equation}{0}  % reset counter 
\setcounter{subsection}{0}
\section*{Appendix B}
%\label{cases}

\subsection{A detailed proof of Theorem \ref{normalize}}

We will assume by induction on $r$  that
\begin{equation}
\g_{\s',\t'} = 1 \text { for } |\t'|=1, |\s'|\leq r-1 \text{ and } \l_{\s}=1
\text{ for } |\s|\leq r
\end{equation}
Fix $\s$ with $|\s|=r+1$.
We need to show that indeed for two decompositions
\begin{equation}
\s = \s'\t'=\s''\t''
\end{equation}
indeed
\begin{equation}
\g_{\s',\t'}= \g_{\s'',\t''}
\end{equation}

We set $\s'''=\s\t'\t''$ and $\t'''=\t'\t''\t'$. It follows
$$
\s'=\s'''\t'', \s'' =\s'''\t''', \t'''\neq \t''
$$

If $|\s'''|=r-1$, we find
\begin{eqnarray*}
\g_{\s',\t'}&=&\g_{\s'''\t'',\t'}\g_{\s''',\t''}
=\g_{\s''',\t''\t'}\g_{\t'',\t'}\\
&=&\g_{\s''',\t'''\t''}\g_{\t''',\t''}
=\g_{\s'''\t''',\t''}\g_{\s''',\t'''}
=\g_{\s'',\t''} 
\end{eqnarray*}

If $|\s'''|=r+1$ then if $\t' =(ij), \quad \t''=(kl)$, 
$i,j,k,l$ must all lie in the same cycle.

Without loss of generality and to avoid too many indices, 
we assume that  this cycle $c$ is just given by 
$
c= (1 2 \cdots h)
$
for some $h\leq r+2$.
First assume that $\{i,j\}\cap \{k,l\}=\emptyset$.
We can then assume $i<j$, $k<l$ and $i<k$.
Then there are three possibilities: 
$i<j<k<l$, $i<k<l<j$ and $i<k<j<l$ where the first two have 
$|\s'''|=r-1$.

So fix $i<k<j<l$.
We see that we can decompose
$$
\s' =\tilde \s (ilh)(kj),\quad  \s''= \tilde \s (ikh)(jl)
$$ 
with
$$
\tilde \s= \s(hljki) \text { and } |\tilde \s|= r-3 
$$

Now
\begin{eqnarray*}
\g_{\s',\t'} &=& \g_{\tilde \s  (ilh)(kj),(ij)}
= \g_{\tilde \s  (ilh)(kj),(ij)} \g_{\tilde \s  (ilh),(kj)}
= \g_{\tilde \s  (ilh),(kj)(ij)} \g_{(kj),(ij)}\\
&=& \g_{\tilde \s  (ilh),(ik)(kj)} \g_{(ik),(kj)}
= \g_{\tilde \s  (ilh)(ik),(kj)} \g_{\tilde \s  (ilh),(ik)}
= \g_{\tilde \s  (iklh),(kj)}\\
&=& \g_{\tilde \s  (ikh)(kl),(kj)} \g_{\tilde \s  (ikh),(kl)}
= \g_{\tilde \s  (ikh),(kl)(kj)} \g_{(kl),(kj)}\\
&=& \g_{\tilde \s  (ikh),(jl)(kl)} \g_{(jl),(kl)}
= \g_{\tilde \s  (ikh)(jl),(kl)} \g_{\tilde \s  (ikh),(jl)}
= \g_{\tilde \s  (ikh)(jl),(kl)}\\
&=& \g_{\s'',\t''}
\end{eqnarray*}
since $|\tilde \s (ilh)|= |\tilde \s (ikh)|= r-1$.

If $|\{i,j\} \cap \{k,l\}|=1$ then we can assume that $j=k$ and 
$i<l$ which leaves us with the cases:
$i<j<l$,$j<i<l$ and $i<l<j$; where in the first two cases 
$|\s'''|= r-1$.

Now assume $i<j<k$.  We can decompose
$$
\s' =\tilde \s (ilh),\quad  \s''= \tilde \s (ijh)
$$ 
with
$$
\tilde \s= \s(hlji) \text { and } |\tilde \s|= h-4 
$$

And
\begin{eqnarray*}
\g_{\s',\t'} &=& \g_{\tilde \s  (ilh),(ij)}
= \g_{\tilde \s  (il)(lh),(ij)} \g_{\tilde \s  (il),(lh)}
= \g_{\tilde \s  (il),(lh)(ij)} \g_{(lh),(ij)}\\
&=& \g_{\tilde \s  (il),(ij)(lh)} \g_{(ij),(lh)}
= \g_{\tilde \s  (il)(ij),(lh)} \g_{\tilde \s  (il),(ij)}
= \g_{\tilde \s  (ijl),(lh)}\nn\\
&=& \g_{\tilde \s  (ij)(jl),(lh)} \g_{\tilde \s  (ij),(jl)}
= \g_{\tilde \s  (ij),(jl)(lh)} \g_{(jl),(lh)}\\
&=& \g_{\tilde \s  (ij),(jh)(jl)} \g_{(jh),(jl)}
= \g_{\tilde \s  (ij)(jh),(jl)} \g_{\tilde \s  (ij),(jh)}
= \g_{\tilde \s  (ijh),(jl)}\\
&=& \g_{\s'',\t''}
\end{eqnarray*}
since $|\tilde \s (il)|= |\tilde \s (ij)|= r-1$.

\end{document}